\documentclass[11pt,draft]{article}
\usepackage[final]{graphicx}            %
\DeclareGraphicsExtensions{.pdf, .jpg}  %
\usepackage{color}

\usepackage{amsmath}                    %
\usepackage{amssymb}                    %
\usepackage{mathtools}                  %
\usepackage{amsthm}
\usepackage{mathrsfs}                   %
\usepackage{stmaryrd}                   %
\usepackage{nth}
\usepackage[symbol,perpage]{footmisc}
\usepackage[sort]{natbib}               %
\usepackage{enumerate}
\usepackage[nosepfour, np, autolanguage]{numprint}
\usepackage{ifthen}
\usepackage{afterpage}
\usepackage[paper = a4paper, centering,
            lmargin = 2.5cm, rmargin = 2.5cm,
            tmargin = 2.5cm, bmargin = 3.0cm]{geometry}
\usepackage{booktabs}
\usepackage{placeins}
\usepackage{shortvrb}

\theoremstyle{plain}                    %
\newtheorem{theorem}{Theorem}           %
\theoremstyle{definition}               %
\newtheorem{problem}{Problem}           %
\theoremstyle{remark}                   %
\usepackage{pdfsync}
\usepackage{subfigure}

\makeatletter

\newif\ifFirstPar       \FirstParfalse

\DeclareRobustCommand\SMC{%
  \ifx\@currsize\normalsize\small\else
   \ifx\@currsize\small\footnotesize\else
    \ifx\@currsize\footnotesize\scriptsize\else
     \ifx\@currsize\large\normalsize\else
      \ifx\@currsize\Large\large\else
       \ifx\@currsize\LARGE\Large\else
        \ifx\@currsize\scriptsize\tiny\else
         \ifx\@currsize\tiny\tiny\else
          \ifx\@currsize\huge\LARGE\else
           \ifx\@currsize\Huge\huge\else
            \small\SMC@unknown@warning
 \fi\fi\fi\fi\fi\fi\fi\fi\fi\fi
}
\newcommand\SMC@unknown@warning{\TBWarning{\string\SMC: unrecognised
    text font size command -- using \string\small}}
\newcommand\textSMC[1]{{\SMC #1}}
\newcommand\acro[1]{\textSMC{#1}\@}

\newcount\T@stCount             \newcount\TestCount
\def\nth#1{%
    \def\reserved@a##1##2\@nil{\ifcat##1n%
           0%
   \let\reserved@b\ensuremath
      \else##1##2%
   \let\reserved@b\relax
      \fi}%
    \TestCount=\reserved@a#1\@nil\relax
    \ifnum\TestCount <0 \multiply\TestCount by\m@ne \fi 
    \T@stCount=\TestCount
    \divide\T@stCount by 100 \multiply\T@stCount by 100
    \advance\TestCount by-\T@stCount     
    \ifnum\TestCount >20 \T@stCount=\TestCount
      \divide\T@stCount by 10 \multiply\T@stCount by 10
      \advance\TestCount by-\T@stCount   
    \fi
     \reserved@b{#1}%
       \textsuperscript{\ifcase\TestCount th
                        \or   st
                        \or   nd
                        \or   rd
                        \else th
                        \fi}%
     }

\makeatother%

\newcommand{\nsd}[1]{\ensuremath{\mathrm{d}{#1}}}
\renewcommand{\d}[1]{\ensuremath{\,\mathrm{d}{#1}}}
\newcommand{\bv}[1]{\boldsymbol{#1}}    %
\newcommand{\tensor}[1]{\mathsf{#1}}    %
\newcommand{\trans}[1]{#1^{\mathrm{T}}} %
\newcommand{\invtrans}[1]{\ensuremath{#1^{-\mathrm{T}}}} %
\newcommand{\FigPath}[1]{figures/#1/#1}
\newcommand{\f}{\text{f}}
\newcommand{\s}{\text{s}}
\DeclareMathOperator{\ldiv}{div}        %
\DeclareMathOperator{\grad}{grad}       %
\DeclareMathOperator{\vssp}{span}        %

\newcommand{\sv}{\ensuremath{1}}
\newcommand{\svs}{\ensuremath{\alpha}}

\newcommand{\sqs}{\ensuremath{\beta}}
\newcommand{\sy}{\ensuremath{3}}
\newcommand{\sys}{\ensuremath{\gamma}}
\newcommand{\map}[1][Blank]{\ifthenelse{\equal{#1}{Blank}}{\ensuremath{\bv{\zeta}(\bv{s},t)}}{\bv{s} + \bv{#1}(\bv{s},t)}}

\title{%
A Fully Coupled Immersed Finite Element Method for Fluid Structure Interaction via the Deal.II Library
}

\author{%
Luca Heltai\footnote{Corresponding Author. Email: Luca Heltai \texttt{<luca.heltai@sissa.it>}; Tel.: +39 040 3787 449; Fax: +39 040 3787528}\\
Scuola Internazionale Superiore di Studi Avanzati\\
                    Via Bonomea 265\\
                    34136 Trieste\\
                    Italy
\and
Saswati Roy \quad and \quad Francesco Costanzo\\
Center for Neural Engineering\\
Department of Engineering Science and Mechanics\\
The Pennsylvania State University\\
University Park
PA 16802
USA}

\begin{document}

\maketitle

\begin{abstract}
We present the implementation of a solution scheme for fluid-structure interaction problems via the finite element software library \texttt{deal.II}.  The solution scheme is an immersed finite element method in which two independent discretizations are used for the fluid and immersed deformable body.  In this type of formulation the support of the equations of motion of the fluid is extended to cover the union of the solid and fluid domains.  The equations of motion over the extended solution domain govern the flow of a fluid under the action of a body force field.  This body force field informs the fluid of the presence of the immersed solid.  The velocity field of the immersed solid is the restriction over the immersed domain of the velocity field in the extended equations of motion.  The focus of this paper is to show how the determination of the motion of the immersed domain is carried out in practice.  We show that our implementation is general, that is, it is not dependent on a specific choice of the finite element spaces over the immersed solid and the extended fluid domains.  We present some preliminary results concerning the accuracy of the proposed method. 
\end{abstract}


\allowdisplaybreaks{                             %


\section{Introduction}
\cite{Heltai2012Variational-Imp0} have recently discussed a fully variational formulation for an immersed method to the solution of fluid-structure interaction (\acro{FSI}) problems. Immersed methods, which deal with the motion of bodies immersed in fluids, allow one to choose the discretization for the fluid and solid domains independently from each other.  As such, they stand in contrast to established methods like the arbitrary Lagrangian-Eulerian (\acro{ALE}) ones (see, e.g., \citealp{HughesLiu_1981_Lagrangian-Eulerian_0}), where the topologies of the solution grids for the fluid and the solid are constrained.

Immersed methods have three main features:
\begin{enumerate}
\item
The support of the equations of motion of the fluid is extended to the union of the physical fluid and solid domains.

\item
The equations of motion of the fluid have terms that, from a continuum mechanics viewpoint, are body forces ``informing'' the fluid of its interaction with the solid.

\item
The velocity field of the immersed solid is identified with the restriction to the solid domain of the velocity field in the equations of motion of the fluid.
\end{enumerate}
A taxonomy of immersed methods can be based on how these three elements are treated theoretically and/or are implemented practically (see the discussion in \citealp{Heltai2012Variational-Imp0}).  Here we employ the approach proposed in \cite{Heltai2012Variational-Imp0} in which the entire solution scheme is developed within the general framework of the finite element method.  Most importantly, the restriction mentioned at point 3 above is done via a fully variational approach.  As such, the approach demonstrated herein stands in contrast to what is used in the immersed boundary methods stemming from the approach of Peskin and his co-workers (see, e.g., \citealp{Peskin_1977_Numerical_0,Peskin_2002_The-immersed_0,GriffithLuo-2012-a,Griffith-2012-a}) or the finite element extension of Peskin's approach due to Liu and co-workers (see, e.g., \citealp{WangLiu_2004_Extended_0, ZhangGerstenberger_2004_Immersed_0,LiuKim_2007_Mathematical_0}), which is based on the implementation of the reproducing kernel particle method.  As explained in detail in \cite{Heltai2012Variational-Imp0}, the method demonstrated herein stems from the approach by \cite{BoffiGastaldi_2003_A-Finite_0}, \cite{Heltai_2006_The-Finite_0}, and \cite{BoffiGastaldiHeltaiPeskin-2008-a}, and~\cite{Heltai-2008-a}.

In Section~\ref{sec: Problem Formulation} we review the problem's governing equations.  In Section~\ref{sec: Reformulation of the governing equations} we present the variational reformulate of the governing equations and we will present their discrete counterparts in Section~\ref{section: discretization by FEM}.   The content of Sections~\ref{sec: Problem Formulation} and~\ref{section: discretization by FEM} follows closely the exposition in \cite{Heltai2012Variational-Imp0} and is reported here for completeness.  In Section~\ref{sec:implementation}  we provide details about the code we have developed and instructions for compilation, execution, and generation of documentation. The entire code is based on the open source \texttt{deal.II} library (see \citealp{BangerthHartmannKanschat-2007-a,BangerthHartmann-deal.II-Differential--0}).  We conclude the article with Section~\ref{sec:numerics}, where  we present some numerical results.

\section{Problem Formulation}
\label{sec: Problem Formulation}

\subsection{Basic notation and governing equations}
\label{subsec: Basic notation and governing equations}
$B_{t}$ in Fig.~\ref{fig: current_configuration}
\begin{figure}[htb]
    \centering
    \includegraphics{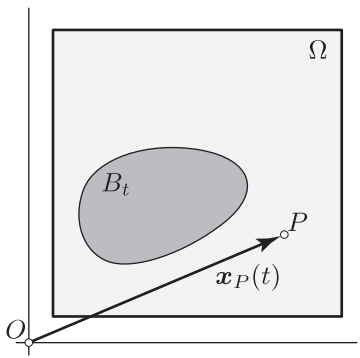}
    \caption{Current configuration $B_{t}$ of a body $\mathscr{B}$ immersed in a fluid occupying the domain $\Omega$.}
    \label{fig: current_configuration}
\end{figure}
represents the configuration of a regular body at time $t$.   $B_{t}$ is a (possibly multiply connected) proper subset of a fixed control volume $\Omega$.  The domain $\Omega\setminus  B_{t}$ is filled by a fluid and we refer to $B_{t}$ as the \emph{immersed body}.  $\partial\Omega$ and $\partial B_{t}$, with outer unit normals $\bv{m}$ and $\bv{n}$, respectively, are the boundaries of $\Omega$ and $B_{t}$. We denote by $B$ the reference configuration of the immersed body.  We denote the position of points of $\mathscr{B}$ in $B$ by $\bv{s}$, whereas we denote the position at time $t$ of a generic point $P \in \Omega$ by $\bv{x}_{P}(t)$.  A motion of $\mathscr{B}$ is a diffeomorphism $\bv{\zeta}: B \to B_{t}$, $\bv{x} = \bv{\zeta}(\bv{s},t)$, with $\bv{s} \in B$, $\bv{x} \in \Omega$, and $t \in [0,T)$, with $T$ a positive real number.

The function $\rho(\bv{x}, t)$ describes the mass density in the entire domain $\Omega$.  The function $\rho$ can be discontinuous across $\partial B_{t}$.  The local form of the balance of mass requires that, $\forall t \in (0,T)$,
\begin{equation}
\label{eq: Balance of mass}
\dot{\rho} + \rho \ldiv\bv{u} = 0,\quad \bv{x} \in \Omega \setminus (\partial\Omega \cup \partial B_{t}),
\end{equation}
where $\bv{u}(\bv{x},t) = \partial\bv{\zeta}(\bv{s},t)/\partial t \big|_{\bv{s} = \bv{\zeta}^{-1}(\bv{x},t)}$ is the velocity field, a dot over a quantity denotes the material time derivative of that quantity, and where `$\ldiv$' represents the divergence operator with respect to $\bv{x}$.

The local form of the momentum balance laws require that, $\forall t \in (0,T)$, $\tensor{T} = \trans{\tensor{T}}$ (the superscript $\mathrm{T}$ denotes the transpose) and
\begin{equation}
\label{eq: Cauchy theorem}
\ldiv \tensor{T} + \rho \bv{b} = \rho \dot{\bv{u}},
\quad \bv{x} \in \Omega \setminus (\partial\Omega \cup \partial B_{t}),
\end{equation}
where $\tensor{T}(\bv{x},t)$ is the Cauchy stress and $\bv{b}(\bv{x},t)$ is the external force density per unit mass acting on the system.

In addition to Eqs.~\eqref{eq: Balance of mass} and~\eqref{eq: Cauchy theorem}, we demand that the velocity field be continuous (corresponding to a no slip condition between solid and fluid) and that the jump condition of the balance of linear momentum be satisfied across $\partial B_{t}$:
\begin{equation}
\label{eq: jump conditions}
\bv{u}(\check{\bv{x}}^{+},t)  = \bv{u}(\check{\bv{x}}^{-},t)
\quad \text{and} \quad
\tensor{T}(\check{\bv{x}}^{+},t) \bv{n} = \tensor{T}(\check{\bv{x}}^{-},t) \bv{n},
\quad \check{\bv{x}} \in \partial B_{t},
\end{equation}
where the superscripts $-$ and $+$ denote limits as $\bv{x} \to \check{\bv{x}}$ from within and without $B_{t}$, respectively.

We denote by $\partial\Omega_{D}$ and $\partial\Omega_{N}$ the subsets of $\partial\Omega$ where Dirichlet and Neumann boundary data are  prescribed, respectively. The domains $\partial\Omega_{D}$ and $\partial\Omega_{N}$ are such that
\begin{equation}
\label{eq: ND boundary}
\partial\Omega = \partial\Omega_{D} \cup \partial\Omega_{N}
\quad \text{and} \quad
\partial\Omega_{D} \cap \partial\Omega_{N} = \emptyset.
\end{equation}
We denote by $\bv{u}_{g}(\bv{x},t)$, with $\bv{x} \in \partial\Omega_{D}$, and by $\bv{\tau}_{g}(\bv{x},t)$, with $\bv{x}\in\partial\Omega_{N}$, the prescribed values of velocity (Dirichlet data) and traction (Neumann data), respectively, i.e.,
\begin{equation}
\label{eq: boundary conditions}
\bv{u}(\bv{x},t) = \bv{u}_{g}(\bv{x},t),\quad \text{for $\bv{x} \in \partial\Omega_{D}$,}
\quad \text{and} \quad
\tensor{T}(\bv{x},t) \bv{m}(\bv{x},t) = \bv{\tau}_{g}(\bv{x},t), \quad \text{for $\bv{x} \in \partial\Omega_{N}$,}
\end{equation}
where the subscript $g$ stands for `given'.

\subsection{Constitutive behavior}
\label{subsec: Constitute behavior}

\paragraph{Constitutive response of the fluid.}
We assume that the fluid is linear viscous and incompressible with uniform mass density $\rho$.  Denoting by $\hat{\tensor{T}}_{\f}$ the constitutive response function of the Cauchy stress of the fluid, we have (see, e.g., \citealp{GurtinFried_2010_The-Mechanics_0})
\begin{equation}
\label{eq: incompressible NS fluid}
\hat{\tensor{T}}_{\f} = -p \tensor{I} + 2 \mu \tensor{D},
\quad
\tensor{D} = \tfrac{1}{2} \bigl(\tensor{L} + \trans{\tensor{L}}\bigr),
\end{equation}
where $p$ is the pressure of the fluid, $\tensor{I}$ is the identity tensor, $\mu > 0$ is a given viscosity coefficient, and $\tensor{L} = \grad \bv{u}$, and where a ``hat'' ($\hat{\tensor{T}}$) is used to distinguish the constitutive response function for $\tensor{T}$ from $\tensor{T}$ itself.  For convenience, we denote by $\hat{\tensor{T}}^{v}_{\f}$ the viscous component of $\hat{\tensor{T}}_{\f}$, i.e.,
\begin{equation}
\label{eq: viscous T fluid}
\hat{\tensor{T}}^{v}_{\f} = 2 \mu \, \tensor{D} = \mu \, \bigl(\tensor{L} + \trans{\tensor{L}}\bigr).
\end{equation}
Due to incompressibility, the balance of mass reduces to
\begin{equation}
\label{eq: incompressibility constraint fluid}
\ldiv \bv{u} = 0 \quad \text{for $\bv{x} \in \Omega\setminus B_{t}$}.
\end{equation}
Under these conditions, $p$ is a Lagrange multiplier that allows to enforce Eq.~\eqref{eq: incompressibility constraint fluid}.  

\paragraph{Constitutive response of the solid.}
The immersed body is taken to be incompressible and viscoelastic of differential type:
\begin{equation}
\label{eq: solid Cauchy Response Function}
\hat{\tensor{T}}_{\s} = -p \tensor{I} +\hat{\tensor{T}}^{e}_{\s} + \hat{\tensor{T}}^{v}_{\s},
\end{equation}
where $\hat{\tensor{T}}^{e}_{\s}$ and $\hat{\tensor{T}}^{v}_{\s}$ denote the elastic and viscous parts of $\hat{\tensor{T}}_{\s}$, respectively, and $p$ is the Lagrange multiplier needed to enforce incompressibility. The viscous part of the behavior is assumed to be of the same type as that of the fluid, that is,
\begin{equation}
\label{eq: viscous part solid}
\hat{\tensor{T}}^{v}_{\s} = 2 \mu \, \tensor{D} = \mu \, \bigl(\tensor{L} + \trans{\tensor{L}} \bigr),
\end{equation}
where $\mu$ is the same constant viscosity coefficient of the fluid.  We assumed that $\hat{\tensor{T}}^{e}_{\s}$ derived from a strain energy potential.  To be precise, let the first Piola-Kirchhoff stress tensor be $\tensor{P}$.  This tensor is related to $\tensor{T}$ as follows (see, e.g., \citealp{GurtinFried_2010_The-Mechanics_0}):
\begin{equation}
\label{eq: P defs}
\tensor{P} = J \tensor{T} \invtrans{\tensor{F}},
\end{equation}
where $J = \det\tensor{F}$, and the tensor $\tensor{F}$, called the deformation gradient, is defined as
\begin{equation}
\label{eq: F defs}
\tensor{F} = \frac{\partial \map}{\partial \bv{s}}.
\end{equation}
Letting $\hat{\tensor{P}}^{e}_{\s} = J \hat{\tensor{T}}^{e}_{\s}\invtrans{F}$ denote the constitutive response function for the elastic  part of the first Piola-Kirchhoff stress tensor, as is typical in  elasticity, we assume the existence of a function $\hat{W}^{e}_{\s}(\tensor{F})$ such that
\begin{equation}
\label{eq: Elastic 1stPK stress}
\hat{\tensor{P}}^{e}_{\s} = \frac{\partial\hat{W}^{e}_{\s}(\tensor{F})}{\partial{\tensor{F}}},
\end{equation}
where $\hat{W}^{e}_{\s}$ is the density of the elastic strain energy of the solid per unit volume.  Invariance under changes of observer demands that $\hat{W}^{e}_{\s}$ be a function of an objective strain measure such as $\tensor{C} = \trans{\tensor{F}}\tensor{F}$.  If the solid is isotropic, $\hat{W}^{e}_{\s}$ must be a function of the principal invariants of $\tensor{C}$.

\subsection{Reformulation of the governing equations}
\label{sec: Reformulation of the governing equations}
We now reformulate the governing equations in variational form.  The motion of the solid will be described via the displacement field, denoted by $\bv{w}$ and defined as
\begin{equation}
\label{eq: S disp def}
\bv{w}(\bv{s},t) := \bv{\zeta}(\bv{s},t) - \bv{s},
\quad \bv{s} \in B.
\end{equation}
The displacement gradient relative to the position in $B$ is denoted by $\tensor{H}$:
\begin{equation}
\label{eq: disp grad def}
\tensor{H} := \frac{\partial \bv{w}}{\partial \bv{s}}
\quad \Rightarrow \quad
\tensor{H} = \tensor{F} - \tensor{I}.
\end{equation}
Equation~\eqref{eq: S disp def} implies
\begin{equation}
\label{eq: w u rel}
\dot{\bv{w}}(\bv{s},t) = \bv{u}(\bv{x},t)\big|_{\bv{x} = \bv{\zeta}(\bv{s},t)}.
\end{equation}
The principal unknowns of our fluid-structure interaction problem are then the fields
\begin{equation}
\label{eq: unknowns}
\bv{u}(\bv{x},t), \quad
p(\bv{x},t), \quad \text{and} \quad
\bv{w}(\bv{s},t),
\quad\text{with $\bv{x} \in \Omega$, $\bv{s} \in B$, and $t \in [0,T)$.}
\end{equation}
The functional spaces for the problem are
\begin{gather}
\label{eq: functional space u}
\bv{u} \in \mathscr{V} = H_{D}^{1}(\Omega)^{d} := \Bigl\{ \bv{u} \in L^{2}(\Omega)^{d} \,\big|\, \nabla_{\bv{x}} \bv{u} \in L^{2} (\Omega)^{d \times d},  \bv{u}|_{\partial\Omega_{D}} = \bv{u}_{g} \Bigr\},
\\
\label{eq: functional space p}
p \in \mathscr{Q} := L^{2}(\Omega), \\
\label{eq: functional space w}
\bv{w} \in \mathscr{Y} = H^{1}(B)^{d} := \Bigl\{ \bv{w} \in L^{2}(B)^{d} \,\big|\, \nabla_{\bv{s}} \bv{w} \in L^{2} (B)^{d \times d} \Bigr\},
\end{gather}
where $\nabla_{\bv{x}}$ and $\nabla_{\bv{s}}$ denote the gradient operators relative to $\bv{x}$ and $\bv{s}$, respectively.
Also, referring to Eq.~\eqref{eq: functional space u}, the function space for the test functions for the velocity field is taken to be as follows:
\begin{equation}
\label{eq: space of test functions v}
\mathscr{V}_{0} = H_{0}^{1}(\Omega)^{d} := \Bigl\{ \bv{v} \in L^{2}(\Omega)^{d} \,\big|\, \nabla_{\bv{x}} \bv{v} \in L^{2} (\Omega)^{d \times d},  \bv{v}|_{\partial\Omega_{D}} = \bv{0} \Bigr\}.
\end{equation}
\subsection{Variational restatement of the governing equations}
\label{subsec: Governing equations: incompressible solid}
When the solid is incompressible, the mass density of both the fluid and the solid are constant so that $\dot{\rho} = 0$ (almost) everywhere in $\Omega$.  Then, referring to Eqs.~\eqref{eq: boundary conditions}, Eqs.~\eqref{eq: functional space u}--\eqref{eq: functional space w}, and the constitutive response functions of both the fluid and the solid, the governing equations introduced so far can be expressed in weak form as follows:
\begin{gather}
\label{eq: Bmomentum weak partitioned first}
\begin{multlined}[b]
\int_{\Omega} \rho(\dot{\bv{u}} - \bv{b}) \cdot \bv{v} \d{v}
+
\int_{\Omega} \hat{\tensor{T}}_{\f} \cdot \nabla_{\bv{x}}\bv{v} \d{v} \\
+
\int_{B_{t}} \bigr(\hat{\tensor{T}}_{\s} - \hat{\tensor{T}}_{\f}\bigl)\cdot \nabla_{\bv{x}}\bv{v} \d{v} - \int_{\partial\Omega_{N}} \bv{\tau}_{g} \cdot \bv{v} \d{a} = 0
\quad \forall \bv{v} \in \mathscr{V}_{0}
\end{multlined}
\shortintertext{and}
\label{eq: Bmass weak partitioned}
\int_{\Omega} q \ldiv \bv{u} \d{v} = 0
\quad \forall q \in \mathscr{Q}.
\end{gather}
A crucial aspect of our approach is the enforcement of Eq.~\eqref{eq: w u rel}.  We enforce this relation weakly as follows:
\begin{equation}
\label{eq: w u rel weak}
\Phi_{B}
\int_{B} \Bigl[\dot{\bv{w}}(\bv{s},t) - \bv{u}(\bv{x},t)\big|_{\bv{x} = \map}\Bigr] \cdot \bv{y}(\bv{s}) \d{V} = 0
\quad
\forall \bv{y} \in \mathscr{Y},
\end{equation}
where $\nsd{V}$ is an infinitesimal volume element of $B$, and where $\Phi_{B}$ is a constant with dimensions of mass over time divided by length cubed, i.e., dimensions such that, in 3D, the volume integral of the quantity $\Phi_{B} \dot{\bv{w}}$ has the same dimensions as a force.  We observe that, since we have assumed that the viscous part of the stress response of the solid is the same as that of the fluid (\citealp{Heltai2012Variational-Imp0} discuss the most general of cases in which the immersed body and the surrounding fluid can have different constitutive response functions), the term $\bigr(\hat{\tensor{T}}_{\s} - \hat{\tensor{T}}_{\f}\bigl)$ in Eq.~(\ref{eq: Bmomentum weak partitioned first}) is equal to the elastic response of the solid $\hat{\tensor{T}}^e_{\s}$. 

Our numerical approximation scheme for  Eqs.~\eqref{eq: Bmomentum weak partitioned first}--\eqref{eq: w u rel weak} is based on the use of two independent triangulations, namely, one of $\Omega$ and one of $B$.  The fields $\bv{u}$ and $p$, as well as their corresponding test functions, will be expressed via finite element spaces supported by the triangulation of $\Omega$.  By contrast, the field $\bv{w}$ will be expressed via a finite element space supported by the triangulation of $B$.  Because of this, any term in Eq.~\eqref{eq: Bmomentum weak partitioned first} defined over $B_{t}$ is now rewritten as an integral over $B$:
\begin{multline}
\label{eq: Bmomentum weak partitioned last}
\int_{\Omega} \rho (\dot{\bv{u}} - \bv{b}) \cdot \bv{v} \d{v}
- \int_{\Omega} p \ldiv \bv{v} \d{v}
+
\int_{\Omega} \hat{\tensor{T}}^{v}_{\f} \cdot \nabla_{\bv{x}}\bv{v} \d{v}
-\int_{\partial\Omega_{N}} \bv{\tau}_{g} \cdot \bv{v} \d{a}
\\
+
\int_{B} \hat{\tensor{P}}^{e}_{\s} \, \trans{\tensor{F}}(\bv{s},t) \cdot \nabla_{\bv{x}}\bv{v}(\bv{x})\bigr|_{\bv{x} = \map} \d{V}
=
0
\quad \forall \bv{v} \in \mathscr{V}_{0}.
\end{multline}

We now define the operators we will use in our finite element formulation.  In these definitions, we will use the following notation:
\begin{equation}
\label{eq: duality notation}
\prescript{}{V^{*}}{\bigl\langle} \psi, \phi \big\rangle_{V},
\end{equation}
in which, given a vector space $V$ and its dual $V^{*}$, $\psi$ and $\phi$ are elements of the vector spaces $V^{*}$ and $V$, respectively, and where $\prescript{}{V^{*}}{\bigl\langle} \bullet, \bullet \big\rangle_{V}$ identifies the duality product between $V^{*}$ and  $V$.  For convenience, we also introduce the following shorthand notation
\begin{align}
\label{eq: fv stress abbreviated}
\hat{\tensor{T}}^{v}[\bv{u}] &= \mu \bigl[\nabla_{\bv{x}}\bv{u}(\bv{x},t) + \trans{(\nabla_{\bv{x}}\bv{u}(\bv{x},t))} \bigr],
\\
\label{eq: Fw abbreviated}
\tensor{F}[\bv{w}] &= \tensor{I} + \nabla_{s}\bv{w}(\bv{s},t),
\\
\label{eq: se stress abbreviated}
\hat{\tensor{P}}^{e}_{\s}[\bv{w}] &= \frac{\partial\hat{W}^{e}_{\s}(\tensor{F})}{\partial{F}}\bigg|_{\tensor{F} = \tensor{F}[\bv{w}]}.
\end{align}
Finally, to help identify the domain and range of these operators, we establish the following convention.  We will use the numbers $1$, $2$, and $3$ to identify the spaces $\mathscr{V}$, $\mathscr{Q}$, and $\mathscr{Y}$, respectively.  We will use the Greek letter $\alpha$, $\beta$, and $\gamma$ to identify the spaces $\mathscr{V}^{*}$, $\mathscr{Q}^{*}$, and $\mathscr{Y}^{*}$, respectively.  Then, a Greek letter followed by a number will identify an operator whose domain is the space corresponding to the number, and whose co-domain is in the space corresponding to the Greek letter.  For example, the notations 
\begin{equation}
\label{eq: space convention}
\mathcal{E}_{\alpha 2}
\quad \text{and} \quad
\mathcal{E}_{\alpha 2} \, p
\end{equation}
will identify a map ($\mathcal{E}_{\alpha 2}$) from $\mathscr{Q}$ into  $\mathscr{V}^{*}$ and the action of this map ($\mathcal{E}_{\alpha 2}\, p \in \mathscr{V}^{*}$) on the field $p \in \mathscr{Q}$, respectively.  If an operators has only one subscript, that subscript identifies the space containing the range of the operator.  With this in mind, let
\begin{alignat}{3}
\label{eq: MOmega def}
\mathcal{M}_{\svs\sv} &: \mathscr{V} \to \mathscr{V}^{*},
&\quad
\prescript{}{\mathscr{V}^{*}}{\bigl\langle}
\mathcal{M}_{\svs\sv}\bv{u},\bv{v}
\big\rangle_{\mathscr{V}} 
&:= 
\int_{\Omega} \rho \, \bv{u} \cdot \bv{v} \d{v}
&\quad
&\forall \bv{u} \in \mathscr{V}, \forall \bv{v} \in \mathscr{V}_{0},
\\
\label{eq: NOmega def}
\mathcal{N}_{\svs\sv}(\bv{u})
&:
\mathscr{V} \to \mathscr{V}^{*},
&\quad
\prescript{}{\mathscr{V}^{*}}{\bigl\langle}
\mathcal{N}_{\svs\sv}(\bv{u}) \bv{w} , \bv{v}
\big\rangle_{\mathscr{V}} 
&:=
\int_{\Omega} \rho (\nabla_{\bv{x}} \bv{w})\bv{u} \cdot \bv{v} \d{v}
&\quad
&\forall \bv{u},\bv{w} \in \mathscr{V}, \forall \bv{v} \in \mathscr{V}_{0},
\\
\label{eq: AOmega def}
\mathcal{D}_{\svs\sv} &: \mathscr{V} \to \mathscr{V}^{*},
&\quad
\prescript{}{\mathscr{V}^{*}}{\bigl\langle}
\mathcal{D}_{\svs\sv}\bv{u},\bv{v}
\big\rangle_{\mathscr{V}} 
&:= 
\int_{\Omega} \hat{\tensor{T}}^{v}_{\s}[\bv{u}] \cdot \nabla_{\bv{x}}\bv{v} \d{v}
&\quad
&\forall \bv{u} \in \mathscr{V}, \forall\bv{v} \in \mathscr{V}_{0},
\\
\mathcal{B}_{\sqs\sv} &: \mathscr{V} \to \mathscr{Q}^{*},
&\quad
\prescript{}{\mathscr{Q}^{*}}{\bigl\langle}
\mathcal{B}_{\sqs\sv} \bv{u}, q
\big\rangle_{\mathscr{Q}} &:= -\int_{\Omega} q \ldiv \bv{u} \d{v}
&\quad
&\forall q \in \mathscr{Q}, \forall \bv{u} \in \mathscr{V}.
\end{alignat}
The operators defined in Eqs.~\eqref{eq: MOmega def}--\eqref{eq: AOmega def} are found in traditional variational formulations of the Navier-Stokes equations and will be referred to as the Navier-Stokes component of the problem.  As typical of other immersed methods, these operators have their support in $\Omega$ as a whole.

We now define the operator in our formulation that has its support over $B$ but does not contain prescribed body forces or boundary terms.
\begin{align}
\label{eq: pseudo stiffness BOmega def}
\begin{split}
&\mathcal{A}_{\svs}(\bv{w},\bv{h}) \in \mathscr{V}^{*},~\forall \bv{w}, \bv{h} \in \mathscr{Y}, \forall \bv{u}\in\mathscr{V}, \forall\bv{v} \in \mathscr{V}_{0}
\\
&\qquad
\begin{aligned}[b]
\prescript{}{\mathscr{V}^{*}}{\bigl\langle}
\mathcal{A}_{\svs}(\bv{w},\bv{h}), \bv{v}
\big\rangle_{\mathscr{V}} &:= 
\int_{B}
\bigl[
\hat{\tensor{P}}^{e}_{\s}[\bv{w}] \trans{\tensor{F}}[\bv{h}]
\cdot \nabla_{\bv{x}} \bv{v}(\bv{x})\bigr]_{\bv{x}=\map[h]} \d{V}.
\end{aligned}
\end{split}
\end{align}

We now define operators with support in $B$ that express the coupling of the velocity fields defined over $\Omega$ and over $B$.  Specifically, we have
\begin{align}
\label{eq: MB def}
\begin{split}
&\mathcal{M}_{\sys\sy} : \mathscr{Y} \to \mathscr{Y}^{*},~\forall \bv{w},\bv{y} \in \mathscr{Y},
\\
&\qquad
\prescript{}{\mathscr{Y}^{*}}{\bigl\langle}
\mathcal{M}_{\sys\sy}\bv{w}, \bv{y}
\big\rangle_{\mathscr{Y}} := \Phi_{B} \int_{B} \bv{w} \cdot \bv{y}(\bv{s}) \d{V},
\end{split}
\\
\label{eq: MGamma def}
\begin{split}
&\mathcal{M}_{\sys\sv}(\bv{w}) : \mathscr{V} \to \mathscr{Y}^{*},~\forall \bv{u} \in \mathscr{V}, \forall \bv{w},\bv{y} \in \mathscr{Y},
\\
&\qquad
\prescript{}{\mathscr{Y}^{*}}{\bigl\langle}
\mathcal{M}_{\sys\sv}(\bv{w}) \bv{u}, \bv{y}
\big\rangle_{\mathscr{Y}} := \Phi_{B} \int_{B} \bv{u}(\bv{x},t)\big|_{\bv{x} = \map[w]} \cdot \bv{y}(\bv{s}) \d{V},
\end{split}
\end{align}

Finally, we define the operators that express the action of prescribed body and surface forces.
\begin{align}
\label{eq: Forcing Omega def}
\begin{split}
&\mathcal{F}_{\svs} \in \mathscr{V}^{*},~\forall \bv{b} \in H^{-1}(\Omega), \forall \bv{\tau}_{g} \in H^{-\frac 1 2}(\partial \Omega_N), \forall \bv{v} \in \mathscr{V}_{0}
\\
&\qquad\prescript{}{\mathscr{V}^{*}}{\bigl\langle}
\mathcal{F}_{\svs}, \bv{v}
\big\rangle_{\mathscr{V}} :=
\int_{\Omega} \rho \, \bv{b} \cdot \bv{v} \d{v} + \int_{\partial\Omega_{N}} \bv{\tau}_g \cdot \bv{v} \d{a}
\end{split}
\\
\label{eq: Forcing B def}
\begin{split}
&\mathcal{G}_{\svs}(\bv{w}) \in \mathscr{V}^{*},~\forall \bv{w} \in \mathscr{Y}, \forall \bv{b} \in H^{-1}(\Omega), \forall \bv{v} \in \mathscr{V}_{0}
\\
&\qquad\prescript{}{\mathscr{V}^{*}}{\bigl\langle}
\mathcal{G}_{\svs}(\bv{w}), \bv{v}
\big\rangle_{\mathscr{V}} :=
\int_{B} \bigl(\rho_{\s_{0}}(\bv{s}) - \rho J[\bv{w}] \bigr) \bv{b} \cdot \bv{v}(\bv{x})\bigr|_{\bv{x} = \map[w]} \d{v}.
\end{split}
\end{align}

In the definition of the operator $\mathcal{A}_{\svs}$ in Eq.~\eqref{eq: pseudo stiffness BOmega def}, the motion of the immersed solid plays a double role in that it affects the elastic response of the solid (through $\bv{w}$) as well as the map (through $\bv{h}$) functioning as a change of variables of integration.  As discussed in \cite{Heltai2012Variational-Imp0}, it is important to separate these two roles and view $\mathcal{A}_{\svs}$ as the composition of a \emph{change of variable} operator and a Lagrangian elastic operator.  To do so, we write
\begin{align}
\label{eq: S def}
\begin{split}
&\mathcal{S}_{\svs\sys}(\bv{h}): \mathscr{H}_{Y}^{*} \to \mathscr{V}^{*},~\forall \bv{y}^{*} \in \mathscr{H}_{Y}^{*}, \forall \bv{h} \in \mathscr{Y}, \forall \bv{v} \in \mathscr{V}_{0}
\\
&\qquad
\prescript{}{\mathscr{V}^{*}}{\bigl\langle}
\mathcal{S}_{\svs\sys}(\bv{h}) \bv{y}^{*},\bv{v}
\big\rangle_{\mathscr{V}} 
:= 
\prescript{}{\mathscr{H}_{Y}^{*}}{\bigl\langle}
\bv{y}^{*},\bv{v}(\bv{x})\big|_{\bv{x} = \bv{s} + \bv{h}(\bv{s})}
\big\rangle_{\mathscr{H}_{Y}},
\end{split}
\\
\label{eq: Agamma def}
\begin{split}
&\mathcal{A}_{\sys}(\bv{w}) \in \mathscr{H}_{Y}^{*},~\forall \bv{w}\in \mathscr{Y}, \forall \bv{y} \in \mathscr{H}_{Y}
\\
&\qquad
\prescript{}{\mathscr{H}_{Y}^{*}}{\bigl\langle}
\mathcal{A}_{\sys}(\bv{w}), \bv{y}
\big\rangle_{\mathscr{H}_{Y}} := \int_{B}
\hat{\tensor{P}}_{\s}^{e}[\bv{w}] \cdot \nabla_{\bv{s}}\bv{y} \d{V}.
\end{split}
\end{align}
Once the operators $\mathcal{S}_{\svs\sys}(\bv{h})$ and $\mathcal{A}_{\sys}(\bv{w})$ are defined, one can prove the following theorem (see \citealp{Heltai2012Variational-Imp0}):
\begin{theorem}[Eulerian and Lagrangian elastic stiffness operators of the immersed domain]
\label{th: eulerian vs lagrangian stiffness}
With reference to the definitions in Eqs.~\eqref{eq: pseudo stiffness BOmega def}, \eqref{eq: S def}, and~\eqref{eq: Agamma def}, we have
\begin{equation}
\label{eq: EulerianLagrangianElasticity}
\mathcal{A}_{\svs}(\bv{w},\bv{h}) = \mathcal{S}_{\svs\sys}(\bv{h}) \mathcal{A}_{\sys}(\bv{w})
\quad \text{and} \quad 
\mathcal{S}_{\svs\sys}(\bv{h}) = \trans{\mathcal{M}}_{\sys\sv}(\bv{h}) \mathcal{M}_{\sys\sy}^{-1},
\end{equation}
where $\mathcal{S}_{\svs\sys}(\bv{h}) \mathcal{A}_{\sys}(\bv{w})$ and $\trans{\mathcal{M}}_{\sys\sv}(\bv{h}) \mathcal{M}_{\sys\sy}^{-1}$ indicate the composition of the operators $\mathcal{S}_{\svs\sys}(\bv{h})$ and $\mathcal{A}_{\sys}(\bv{w})$ and of the operators $\trans{\mathcal{M}}_{\sys\sv}(\bv{h})$ and $\mathcal{M}_{\sys\sy}^{-1}$, respectively.
\end{theorem}
The operators defined above allow us to formally restate the overall problem described by Eqs.~\eqref{eq: Bmomentum weak partitioned last}, \eqref{eq: Bmass weak partitioned}, and~\eqref{eq: w u rel weak} as follows:
\begin{problem}[Dual formulation]
\label{prob: IFIS}
Given initial conditions $\bv{u}_{0} \in \mathscr{V}$ and $\bv{w}_{0} \in \mathscr{Y}$, for all $t \in (0,T)$ find $\bv{u}(\bv{x},t) \in \mathscr{V}$, $p(\bv{x},t) \in \mathscr{Q}$, and $\bv{w}(\bv{s},t) \in \mathscr{Y}$ such that
\begin{align}
\label{eq: BLM Formal dual}
\mathcal{M}_{\svs\sv}\bv{u}' +  \mathcal{N}_{\svs\sv}(\bv{u})\bv{u} + 
\mathcal{D}_{\svs\sv}\bv{u} + \trans{(\mathcal{B}_{\sqs\sv})}p
+ \mathcal{S}_{\svs\sys}(\bv{w}) \mathcal{A}_{\sys}(\bv{w})
&= \mathcal{F}_{\svs} + \mathcal{G}_{\svs}(\bv{w}),
\\
\label{eq: incompressibility Formal dual}
\mathcal{B}_{\sqs\sv}\bv{u} &= 0,
\\
\label{eq: velocity coupling dual}
\mathcal{M}_{\sys\sy}\bv{w}' - \mathcal{M}_{\sys\sv}(\bv{w})\bv{u} &= \bv{0},
\end{align}
where $\bv{u}'(\bv{x},t) = \partial\bv{u}(\bv{x},t)/\partial t$ and $\bv{w}'(\bv{s},t) = \partial \bv{w}(\bv{s},t)/\partial t$.
\end{problem}

Problem~\ref{prob: IFIS} can be formally presented in terms of the Hilbert space $\mathscr{Z} := \mathscr{V}\times \mathscr{Q}\times \mathscr{Y}$, and $\mathscr{Z}_{0} := \mathscr{V}_{0} \times \mathscr{Q} \times \mathscr{Y}$ with inner product given by the sum of the inner products of the generating spaces. Defining $\mathscr{Z} \ni \xi := \trans{[\bv{u}, p, \bv{w}]}$ and $\mathscr{Z}_{0} \ni \psi := \trans{[\bv{v}, q, \bv{y}]}$, then Problem~\ref{prob: IFIS} can be compactly stated as
\begin{problem}[Grouped dual formulation]
\label{prob: IFG}
Given an initial condition $\xi_0 \in \mathscr{Z}$, for all $t \in (0,T)$ find $\xi(t) \in \mathscr{Z}$, such that
\begin{equation}
\label{eq:formal grouped dual}
\langle \mathcal{F}(t, \xi, \xi') , \psi \rangle =0, \quad \forall \psi \in \mathscr{Z}_0,
\end{equation}
where the full expression of $\mathcal{F} : \mathscr{Z} \mapsto \mathscr{Z}_0^*$ is defined as in Problem~\ref{prob: IFIS}.
\end{problem}
The energy estimates concerning the above abstract formulation has been discussed in \cite{Heltai2012Variational-Imp0} where it is shown that stability is obtained under the same assumptions that yield stability for the Navier-Stokes problems.

\section{Discretization}
\label{section: discretization by FEM}
\subsection{Spatial discretization by finite elements}
The fluid domain is discretized into the triangulation $\Omega_{h}$ and the immersed body into the trangulation $B_{h}$ .  each of these triangulations consists of (closed) cells $K$ (triangles or quadrilaterals in 2D, and tetrahedra or hexahedra in 3D) such that:
\begin{enumerate}
\item
$\overline{\Omega} = \cup \{ K \in \Omega_{h} \}$, and $\overline{B} = \cup \{ K \in B_{h} \}$;

\item
Any two cells $K,K'$ only intersect in common faces, edges, or vertices;

\item
The decomposition $\Omega_{h}$ matches the decomposition $\partial \Omega = \partial\Omega_{D} \cup \partial\Omega_{N}$.
\end{enumerate}
On $\Omega_{h}$ and $B_{h}$, we define the finite dimensional subspaces $\mathscr{V}_{h} \subset \mathscr{V}$, $\mathscr{Q}_{h} \subset \mathscr{Q}$, and $\mathscr{Y}_{h} \subset \mathscr{Y}$ as follows:
\begin{alignat}{5}
\label{eq: functional space u h}
\mathscr{V}_h &:= \Bigl\{ \bv{u}_h \in \mathscr{V} \,&&\big|\, \bv{u}_{h|K}  &&\in \mathcal{P}_V(K), \, K &&\in \Omega_h \Bigr\} &&\equiv \vssp\{ \bv{v}_{h}^{i} \}_{i=1}^{N_{V}}
\\
\label{eq: functional space p h}
\mathscr{Q}_h &:= \Bigl\{ p_h \in \mathscr{Q} \,&&\big|\, p_{h|K}  &&\in \mathcal{P}_Q(K), \, K &&\in \Omega_h \Bigr\} &&\equiv \vssp\{ q_{h}^{i} \}_{i=1}^{N_{Q}}\\
\label{eq: functional space w h}
\mathscr{Y}_{h} &:= \Bigl\{ \bv{w}_h \in \mathscr{Y} \,&&\big|\, \bv{w}_{h|K} &&\in \mathcal{P}_Y(K), \, K &&\in B_h \Bigr\} &&\equiv \vssp\{ \bv{y}_{h}^{i} \}_{i=1}^{N_{Y}},
\end{alignat}
where $\mathcal{P}_{V}(K)$, $\mathcal{P}_{Q}(K)$ and $\mathcal{P}_{Y}(K)$ are polynomial spaces of degree $r_{V}$, $r_{Q}$ and $r_{Y}$ respectively on the cells $K$, and $N_V$, $N_Q$ and $N_Y$ are the dimensions of each finite dimensional space.  The pair $\mathscr{V}_{h}$ and $\mathscr{Q}_{h}$ are chosen so that the inf-sup condition for the well-posedness of the Navier-Stokes problem (see, e.g., \citealp{BrezziFortin-1991-a}) is satisfied.

The discrete version of Problem~\ref{prob: IFIS} is now presented using a matrix notation.  An element of a discrete space, say $\bv{u}_{h} \in \mathscr{V}_{h}$, is represented by a column vector of time dependent coefficients $u_{h}^{j}(t)$, $j = 1,\ldots,N_V$, such that $\bv{u}_{h}(\bv{x},t) = \sum u_{h}^{j}(t) \bv{v}_{h}^{j}(\bv{x})$, where $\bv{v}_{h}^{j}$ is the $\nth{j}$ base element of $\mathscr{V}_{h}$.  We use the notation $M_{\svs\sv} \bv{u}_h$ to represent the multiplication of the column vector $\bv{u}_h$ by the matrix whose elements $M_{\svs\sv}^{ij}$ are
\begin{equation}
  \label{eq:eq:notation matrix element}
  M_{\svs\sv}^{ij} := 
  \prescript{}{\mathscr{V}^{*}}{\bigl\langle}
  \mathcal{M}_{\svs\sv}\bv{v}_h^j,\bv{v}_h^i
  \big\rangle_{\mathscr{V}},
\end{equation}
where the operator in angle brackets is the one defined earlier. A similar notation is adopted for all other previously defined operators.
With this notation, the duality products in the discrete spaces are indicated by simple scalar products in $\mathbb{R}^N$ ($N$ depending on the dimension of the system at hand). Hence, using the matrix $M_{\svs\sv}$, we can write
  \begin{equation}
    \label{eq:duality product in discrete spaces}
    \prescript{}{\mathscr{V}^{*}}{\bigl\langle}
    \mathcal{M}_{\svs\sv}\bv{u}_h,\bv{v}_h
    \big\rangle_{\mathscr{V}} = \bv{v}_h \cdot M_{\svs\sv} \bv{u}_h,
  \end{equation}
where the dot-product on the right hand side is the scalar product in $\mathbb{R}^{N_V}$.

Chosen $\Omega_{h}$ and $B_{h}$ along with $\mathscr{V}_{h}$, $\mathscr{Q}_{h}$, and $\mathscr{Y}_{h}$, Problem~\ref{prob: IFIS} is reformulated as follows:
\begin{problem}
\label{prob: IFIS - discrete}
Given $\bv{u}_{0} \in \mathscr{V}_{h}$, $\bv{w}_{0} \in \mathscr{Y}_{h}$, for all $t \in (0,T)$, find $\bv{u}_{h}(t) \in \mathscr{V}_{h}$, $p_{h}(t) \in \mathscr{Q}_{h}$, and $\bv{w}_{h}(t) \in \mathscr{Y}_{h}$ such that
\begin{align}
\label{eq: BLM Formal dual discrete}
\begin{multlined}[b][9cm]
{M}_{\svs\sv}\bv{u}_{h}' +  {N}_{\svs\sv}(\bv{u}_{h})\bv{u}_{h} + 
{D}_{\svs\sv}\bv{u}_{h} + \trans{({B}_{\sqs\sv})}p_{h}
+  S_{\svs\sys}(\bv{w}_{h}){A}_{\sys}(\bv{w}_{h})
\end{multlined}
&= {F}_{\svs} + {G}_{\svs}(\bv{w}_{h}),
\\
\label{eq: incompressibility Formal dual discrete}
{B}_{\sqs\sv}\bv{u}_{h} &= 0,
\\
\label{eq: velocity coupling dual discrete}
{M}_{\sys\sy}\bv{w}_{h}' - {M}_{\sys\sv}(\bv{w}_{h})\bv{u}_{h} &= \bv{0},
\end{align}
where $\bv{u}'_{h}(\bv{x},t) = \sum [u_h^j(t)]' \bv{v}_h^j(\bv{x})$ and $\bv{w}'_{h}(\bv{s},t) = \sum [w_h^j(t)]' \bv{y}_h^j(\bv{s})$, and where the prime denotes ordinary differentiation with respect to time.
\end{problem}

In compact notation, Problem~\ref{prob: IFIS - discrete} can be casted as semi-discrete problems in the space $\mathscr{Z} \supset \mathscr{Z}_h := \mathscr{V}_h \times \mathscr{Q}_h \times \mathscr{Y}_h$ as 
\begin{problem}
\label{prob: Combo prob discrete}
Given an initial condition $\xi_0 \in \mathscr{Z}_h$, for all $t \in (0,T)$ find $\xi_h(t) \in \mathscr{Z}_h$, such that
\begin{equation}
  \label{eq: dae formulation}
  F(t, \xi_h, \xi_h') = 0,
\end{equation}
where
\begin{equation}
  \label{eq: dae formulation bis}
  F^i(t, \xi_h, \xi_h') := \langle \mathcal{F}(t, \xi_h, \xi_h') , \psi^i_h \rangle, \quad i=0, \dots, N_V+N_Q+N_Y,
\end{equation}
and $\mathcal{F}$ has the same meaning as in Eq.~\eqref{eq:formal grouped dual}, with $\psi^i_h$ being the basis function for the spaces $\mathscr{V}_h$, $\mathscr{Q}_h$, or $\mathscr{Y}_h$ corresponding to the given value of $i$.
\end{problem}

\subsection{Coupling of the fluid and immersed domains}
\label{subsection: FEM implementation}
The operators $M_{\svs\sv}$, $N_{\svs\sv}(\bv{u}_{h})$, $D_{\svs\sv}$, $B_{\sqs\sv}$, and $F_{\svs}$ in Problem~\ref{prob: IFIS - discrete} are common in variational formulations of the Navier-Stokes problem and were implemented in a standard fashion.  The operator $M_{\sys\sy}$ was also implemented in a standard fashion since it is the mass matrix for $\mathscr{Y}_{h}$.  Less common are the operators that depend nonlinearly on the motion of the immersed domain $\bv{w}$.  Thus, we now discuss the practical implementation of such operators.

Let's consider, for example, the matrix $M_{\sys\sv}(\bv{w})$ contributing to the velocity coupling between the fluid and immersed domain:
\begin{equation}
  \label{eq: MGamma def bis}
  {M}_{\sys\sv}^{ij}(\bv{w}_{h}) 
=
  \prescript{}{\mathscr{Y}^{*}_{h}}{\bigl\langle}
  \mathcal{M}_{\sys\sv}(\bv{w}_{h}) \bv{v}^j_{h}, \bv{y}^i_{h}
  \big\rangle_{\mathscr{Y}_{h}}
= \Phi_{B} \int_{B}  \bv{v}^j_{h}(\bv{x})\big|_{\bv{x} = \bv{s} + \bv{w}_{h}(\bv{s},t)} \cdot \bv{y}^i_{h}(\bv{s}) \d{V}.
%
\end{equation}
The above integral is computed by summing contributions from each cell $K$ of $B_{h}$.  Each of these contributions is a sum over the $N_Q$ quadrature points.  We observe that the integrand $\bv{y}^i_{h}(\bv{s})$ is supported over the triangulation of $B_{h}$ but the functions $\bv{v}^j_{h}(\bv{x})$ (with $\bv{x} = \bv{s} + \bv{w}_{h}(\bv{s},t)$) are supported over the triangulation $\Omega_{h}$.  Therefore, the construction of operators like ${M}_{\sys\sv}^{ij}(\bv{w}_{h})$ draws information from two independent triangulations.
\begin{figure}[htb]
    \centering
    \includegraphics{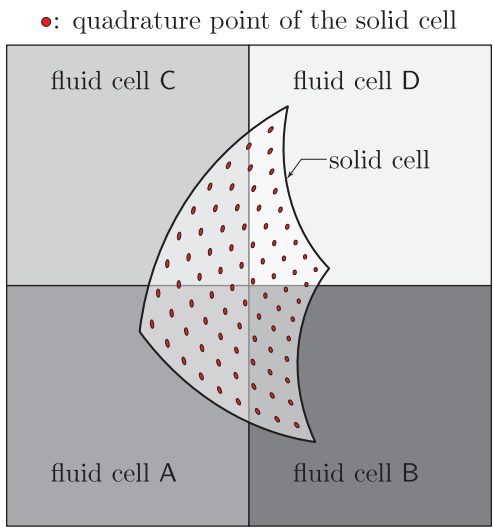}
    \caption{Cells denote as \textsf{A}--\textsf{D} represent a four-cell patch of the triangulation of the fluid domain.  The cell denoted as ``solid cell'' represents a cell of the triangulation of the immersed solid domain that is contained in the union of cells \textsf{A}--\textsf{D} of the fluid domain.  The filled dots represent the quadrature points of the quadrature rule adopted to carry out integration over the cells of the immersed domain.}
    \label{fig: integration}
\end{figure}
In our code, we start by determining the position of the quadrature points of the immersed element, both relative to the reference unit element and relative to the global coordinate system adopted for the calculation, through the mappings:
\begin{alignat}{3}
  \label{eq:mapping Khat K solid}
  \bv{s}_K & : \hat{K} := [0,1]^d &&\mapsto K \in B_h, \\
  \label{eq:mapping K K solid}
  I+\bv{w}_h & : K && \mapsto \text{solid cell}.
\end{alignat}
These maps allow us to determine the global coordinates of the quadrature points. These coordinates are then passed to a search algorithm that identifies the cells in $\Omega_{h}$ that contain the points in question.  In turn, this identification allows is to evaluate the functions $\bv{v}_h^j$. The overall operation is illustrated in Fig.~\ref{fig: integration} where we show a cell of $B_{h}$ straddling four cells of $\Omega_{h}$ denoted fluid cells \textsf{A}--\textsf{D}.  The quadrature points over the solid cell are denoted by filled circles.  The contribution to the integral in Eq.~\eqref{eq: MGamma def bis} due to the solid cell is then computed by summing the partial contributions corresponding to each of the fluid cells intersecting the solid cell in question:
\begin{equation}
  \begin{aligned}[b]
  \label{eq: MGamma def tris}
  {M}_{\sys\sv}^{ij}(\bv{w}_{h}) & =
  \sum_{K\in B_h} \int_{K} \bv{v}^j_{h}(\bv{x})\big|_{\bv{x} = \bv{s} + \bv{w}_{h}(\bv{s},t)} \cdot \bv{y}^i_{h}(\bv{s}) \d{V},
  \\
  & \sim \sum_{K\in B_h} \sum_{q=1}^{N_{K,q}} \bv{v}^j_{h}(\bv{x})\big|_{\bv{x} = \bv{s_{K,q}} + \bv{w}_{h}(\bv{s_{K,q}},t)} \cdot \bv{y}^i_{h}(\bv{s}_{K,q}) \omega_{K,q},
 \end{aligned}
\end{equation}
where $\bv{s}_{K,q}$ is the image of $q$-th quadrature point under the mapping $\bv{s}_K$, and $\omega_{K,q}$ is the corresponding quadrature weight.  The implementation of an efficient search algorithm responsible for identifying the fluid cells intersecting an individual solid cell is the only technically challenging part of the procedure.  We use the built-in facilities of the \texttt{deal.II} library to perform this task. Once the fluid cells containing the quadrature points of a given solid cell are found, we determine the value of $\bv{v}^j_{h}$ at the quadrature points using the interpolation infrastructure inherent in the finite element representation of fields defined over $\Omega_{h}$.  The \texttt{deal.II} C++ class we use for this implementation is the \texttt{FEFieldFunction}.

\subsection{Time discretization}
\label{sec:time-discretization}
Equation~\eqref{eq: dae formulation} represents a system of nonlinear differential algebraic equations (\acro{DAE}), which we solve using a Newton iteration. In the code accompanying this paper, the time derivative $\xi'$ is approximated very simply via an implicit-Euler scheme:
  \begin{equation}
    \label{eq:time derivative theta}
    \xi_n' = h^{-1} \bigl( \xi_{n} - \xi_{n-1} \bigr),
  \end{equation}
where $\xi_{n}$ and $\xi_{n}'$ are the computed approximations to $\xi(t_{n})$ and $\xi'(t_{n})$, respectively, and the step size $h = t_{n} - t_{n-1}$ is kept constant throughout the computation.  Although not second order accurate, this time stepping scheme is asymptotically stable.

The application of the implicit-Euler scheme in Eq.~\eqref{eq:time derivative theta} to the \acro{DAE} system in Eq.~\eqref{eq: dae formulation} results in a nonlinear algebraic system to be solved at each step:
  \begin{equation}
    \label{eq:dae algebraic system}
    G(\xi_{n}) := F\biggl(t_{n}, \xi_{n},  h^{-1} \big( \xi_{n} - \xi_{n - 1} \big)\biggr) = 0.
  \end{equation}
The nonlinear system in Eq.~\eqref{eq:dae algebraic system} is solved via Newton iterations. This leads to a linear system for each Newton correction, of the form
  \begin{equation}
    \label{eq:dae newton correction}
    J[\xi_{n,m+1}-\xi_{n,m}] = -G(\xi_{n,m}),
  \end{equation}
where $\xi_{n,m}$ is the $m$th approximation to $\xi_{m}$. Here $J$ is  some approximation to the system's Jacobian
\begin{equation}
  \label{eq:dae Jacobian}
  J = \frac{\partial G}{\partial \xi} = \frac{\partial F}{\partial \xi} +\alpha \frac{\partial F}{\partial \xi'},
\end{equation}
where $\alpha = 1/h$.
In our finite element implementation, we assemble the residual $G(\xi_{n,m})$ at each Newton correction.  The implementation of the residual vector is based on the formulation presented in Problem~\ref{prob: IFIS - discrete}.  However, this formulation makes the determination of the corresponding Jacobian rather involved due to the structure of the operator $\mathcal{S}_{\svs\sys}(\bv{w})$ (see Eq.~\eqref{eq: EulerianLagrangianElasticity}).  Hence, we have implemented a Newton-Raphson iteration based on an approximate Jacobian.  With reference to Theorem~\ref{th: eulerian vs lagrangian stiffness} and Eq.~\eqref{eq: BLM Formal dual}, the Jacobian we assemble is the exact Jacobian of a formulation in which the operator product $\mathcal{S}_{\svs\sys}(\bv{w}) \mathcal{A}_{\sys}(\bv{w})$ is replaced by the operator $\mathcal{A}_{\svs}(\bv{w},\bv{h})$ defined in Eq.~\eqref{eq: pseudo stiffness BOmega def}.  In the code accompanying this paper, the final system is solved using the direct solver provided by the \texttt{UMFPACK} package (see \citealp{Davis_2004_UMFPACK}).

\section{Implementation}
\label{sec:implementation}
\MakeShortVerb{\|}
\subsection{Source files and library requirements}
\label{sec:source-files-library}

The included source code is based on the |deal.II 7.1.0| library (see \citealp{BangerthHartmann-deal.II-Differential--0}). In what follows, we assume that the user has installed the |deal.II| library in some directory (we have tested our code with |deal.II| version |7.1| and with the latest svn version |7.2 pre|), and all paths will be relative to the install directory path, which we call |deal.II|.  For the program to work properly, |deal.II| should be configured at least with the flag |--with-umfpack|, to enable UMFPACK solver developed in \cite{Davis_2004_UMFPACK} and used extensively inside the program.
  
The provided zip archive should be unzipped under \\
|deal.II/examples/step-feibm|

\begin{table}
\caption{\label{tab:sourcefiles}Content of the provided \texttt{zip} archive.}
\centering
\begin{tabular}{ll}
\toprule
  |INSTALL.txt|: & installation instructions;\\
  |Makefile|: & standard |deal.II| makefile;\\
  |step-feibm.cc|:& main program file;\\
  |immersed_fem.prm|: &default parameter file;\\
  |meshes/|: &directory containing a collection of input mesh files, in
  |UCD| format;\\
  |prms/|: &examples parameter files;\\
|out/:| &empty directory, in which output files will be
  written;\\
|doc/|: &|Doxygen| documentation directory.\\
\bottomrule
\end{tabular}
\end{table}
Table~\ref{tab:sourcefiles} provides a summary of the distributed files and directories.  If the code is unzipped in the above location, it can be compiled by simply typing |make| at the
command line prompt, and run with \\
|./step-feibm parameter_file.prm|

If the file |parameter_file.prm| does not exist, the program creates one with default values, which can then be suitably modified by the user should edit for her needs. We distribute all the parameter files that were used to produce the results in Section~\ref{sec:numerics} along with the needed mesh files. These can be found in the directories |prms| and |meshes|, respectively.
If the program is run without arguments, it is assumed that the problem parameters are those in the file |parameter_file.prm|.  As mentioned earlier, if the file in question does not already exist, a default copy will be created.

If the user has |Doxygen|, a complete and browsable documentation of the source code itself can be generated in two ways: by typing |make| in the subdirectory\\
|deal.II/examples/step-feibm/doc|,\\ 
or by typing |make online-doc| in the directory |deal.II|.  In the first case, the documentation will be accessible at the address\\
|deal.II/examples/step-feibm/doc/html/index.html|,\\
while in the second case, the command will generate the full |deal.II| documentation, together with the documentation of the |step-feibm| source code, which will be
accessible at the address \\
|deal.II/doc/doxygen/deal.II/step_feibm.html|.\\
In the second case, the process is a lot longer but the documentation generated thus has the added benefits of being fully integrated with the |deal.II| documentation as well as having hyperlinks to all the |deal.II| classes that have been used in our program.
\subsection{Parameter and input files}
\label{sec:param-input-files}

The behavior of the program is controlled by the |ProblemParameters<dim>| class, which is derived from the |deal.II| class |ParameterHandler| and is used to define and to read from a file all the problem parameters that the user can set.

The following is a sample parameter file that can be used with our code.

{\footnotesize%
\begin{verbatim}
# Listing of Parameters
# ---------------------

# Time Stepping 
set Final t                                 = 1
set Delta t                                 = .1
set Interval (of time-steps) between output = 1

# Non linear solver 
set Force J update at step beginning        = false
set Update J cont                           = false
set Semi-implicit scheme                    = true
set Use spread operator                     = true

# Constitutive models available are: INH_0: incompressible Neo-Hookean with
# P^{e} = mu (F - F^{-T}); INH_1: incompressible neo-Hookean with P^{e} = mu F; 
# CircumferentialFiberModel: incompressible with P^{e} = mu F
# (e_{\theta} \otimes e_{\theta}) F^{-T}; this is suitable for annular solid
# comprising inextensible circumferential fibers
set Solid constitutive model                = INH_0
set Density                                 = 1
set Viscosity                               = 1
set Elastic modulus                         = 1
# Dimensional constant for the velocity equation
set Phi_B                                   = 1

# Solid mesh information
set Solid mesh                              = meshes/solid_square.inp
set Solid refinement                        = 1

# Fluid mesh information 
set Fluid mesh                              = meshes/fluid_square.inp
set Fluid refinement                        = 4
set All Dirichlet BC                        = true
set Dirichlet BC indicator                  = 1
set Velocity finite element degree          = 2
# Select between FE_Q (Lagrange finite element space of continuous, piecewise
# polynomials) or FE_DGP(Discontinuous finite elements based on Legendre
# polynomials) to approximate the pressure field
set Finite element for pressure             = FE_DGP
set Fix one dof of p                        = false

# Base name used for the output files
set Output base name                        = out/square

# This section is used only when the constitutive model is set to
# CircumferentialFiberModel
subsection Equilibrium Solution of Ring with Circumferential Fibers
  set Any edge length of the (square) control volume = 1.
  set Inner radius of the ring                       = 0.25
  set Width of the ring                              = 0.0625
  set x-coordinate of the center of the ring         = 0.5
  set y-coordinate of the center of the ring         = 0.5
end


subsection W0
  # Sometimes it is convenient to use symbolic constants in the expression
  # that describes the function, rather than having to use its numeric value
  # everywhere the constant appears. These values can be defined using this
  # parameter, in the form `var1=value1, var2=value2, ...'.
  # 
  # A typical example would be to set this runtime parameter to
  # `pi=3.1415926536' and then use `pi' in the expression of the actual
  # formula. (That said, for convenience this class actually defines both `pi'
  # and `Pi' by default, but you get the idea.)
  set Function constants  = 

  # The formula that denotes the function you want to evaluate for particular
  # values of the independent variables. This expression may contain any of
  # the usual operations such as addition or multiplication, as well as all of
  # the common functions such as `sin' or `cos'. In addition, it may contain
  # expressions like `if(x>0, 1, -1)' where the expression evaluates to the
  # second argument if the first argument is true, and to the third argument
  # otherwise. For a full overview of possible expressions accepted see the
  # documentation of the fparser library.
  # 
  # If the function you are describing represents a vector-valued function
  # with multiple components, then separate the expressions for individual
  # components by a semicolon.
  set Function expression = 0; 0

  # The name of the variables as they will be used in the function, separated
  # by commas. By default, the names of variables at which the function will
  # be evaluated is `x' (in 1d), `x,y' (in 2d) or `x,y,z' (in 3d) for spatial
  # coordinates and `t' for time. You can then use these variable names in
  # your function expression and they will be replaced by the values of these
  # variables at which the function is currently evaluated. However, you can
  # also choose a different set of names for the independent variables at
  # which to evaluate your function expression. For example, if you work in
  # spherical coordinates, you may wish to set this input parameter to
  # `r,phi,theta,t' and then use these variable names in your function
  # expression.
  set Variable names      = x,y,t
end


subsection force
  set Function constants  = 
  set Function expression = 0; 0; 0
  set Variable names      = x,y,t
end


subsection u0
  set Function constants  = 
  set Function expression = 0; 0; 0
  set Variable names      = x,y,t
end


subsection ug
  set Function constants  = 
  set Function expression = if(y>.99, 1, 0); 0; 0
  set Variable names      = x,y,t
end
\end{verbatim}}
At the beginning of the parameter we find specifications for the time stepper and for the nonlinear solver.  In addition, we find information on the constitutive behavior of both the fluid and the immersed solid.

The user can specify the names of the files containing the meshes for the control volume and the immersed solid, along with the initial global refinement level for each mesh, in the parameter file.  In the section pertaining to the control volume, the user can also set the degree of the finite element spaces for the fluid velocity as well as the type of the finite element space for the fluid pressure. The type and degree of the finite element space for the displacement of the immersed domain are automatically set to be the same as those for the velocity of the fluid. A degree greater than or equal to two should be selected for the finite element space of the velocity so as to ensure proper inf-sup stability. The degree of the pressure space is then automatically set to be one less than that for the velocity.

In the latter part of the parameter file, the user can specify the initial and boundary values of the solution as well as the external body forces. Here |W0| denotes the initial value of the displacement of the immersed domain, |force| denotes the external body force field, |u0| is the initial condition for the velocity and the pressure fields and |ug| is the Dirichlet boundary condition (here configured for a \emph{lid-cavity} problem).

The above file, for example, generates the parameters for a \emph{lid-cavity} problem inside a square control volume (read from |meshes/fluid_square.inp|), with an immersed solid whose mesh is given in\\
|meshes/solid_square.inp|.

\subsection{Code structure}
\label{sec:code-structure}

The structure of our program follows closely the structure of most tutorial programs in the |deal.II| library, to which we refer for further explanations and examples. The main class of the program is the class |ImmersedFEM<dim>|, in which all objects and methods to solve the problem at hand are defined (including an object of type |ProblemParameters<dim>|).

Execution of the solution is triggered in the method |run()|, which starts the time stepping scheme of the DAE system described in Section~\ref{sec:time-discretization}, and controls the convergence of the Newton iteration scheme for the solution of system Eq.~(\ref{eq:dae newton correction}).

Detailed documentation of the code has been embodied in the code itself, and can be automatically generated with |Doxygen|. Here we only briefly overview the main ideas behind the use of |deal.II| for immersed methods.

Due to the nature of the method, two different sets of objects are needed to describe the triangulation, the degrees of freedom, etc., of both the fluid and the immersed domains. In the code, objects pertaining to the fluid have been denoted with the suffix |_f|, whereas objects pertaining to the immersed solid have been denoted with the suffix |_s|. For example, |tria_s| and |tria_f| are the two |Triangulation<dim>| objects of the solid and fluid domains, respectively.

In the code, solution vectors and residuals are constructed as\\ |BlockVector<double>| objects and the Jacobian matrix is constructed as a \\ |BlockSparseMatrix<double>| object.  This has been done to reflect the logical splitting of these entities between the fluid and the solid, and to allow access to the individual blocks at the same time. We split the vectors and matrices into two and four parts, respectively.  The block vectors storing the overall solutions at the current time step and at the previous time step are called |xi| and |previous_xi|, respectively.  The first block of these block vectors pertains to the fluid and it is of size |n_dofs_up|, which is also equal to |dh_f.n_dofs()|.  The second block pertains to the solid and has a size of |n_dofs_W|, which is also equal to |dh_s.n_dofs()|. 

The various tutorial examples of the |deal.II| library describe in an exhaustive manner how to treat a single triangulation and a single degrees-of-freedom handler for both fluid-only problems (e.g., the example program |step-35|) and elasticity-only problems (e.g., the example program |step-44|). The most delicate part of immersed methods, however, requires the coupling between a fixed background mesh (the fluid), and a moving and deforming foreground mesh (the elastic solid). The deformation of the foreground mesh is achieved very effectively through the \\ |MappingQEulerian<dim,spacedim>| class, which uses the information stored in the displacement vector to automatically compute the deformed positions of the mesh and of the quadrature points in a Lagrangian way. Notice that while the name suggests an Eulerian description, this object in reality performs a Lagrangian iso-parametric transformation from the reference grid, stored in |tria_s|, to the current configuration of the solid via the deformation vector $\bv{w}$. Details on construction and use of this class are given in Section~\ref{sec:immersed-map}.

Evaluation of the quadrature points of the solid on the background fluid mesh is achieved through the class |FEFieldFunction<dim>|, which allows one to evaluate the values of finite element fields at \emph{arbitrary} points. In particular, its method\\
|FEFieldFunction<dim>::compute_point_locations| is the one that returns the lists required to compute the coupling integrals (see Section~\ref{subsection: FEM implementation}) and is used both in the creation of the sparsity pattern that features the coupling between the degrees of freedom of the fluid and the immersed solid (see Section~\ref{sec:sparsity-pattern}), as well as in the assembling of the residual vector and the Jacobian matrix of the DAE system (see Section~\ref{sec:residual-jacobian}).

\subsubsection{Immersed map}
\label{sec:immersed-map}

Whenever it is necessary to compute the deformed configuration of the solid, an iso-parametric displacement is superimposed on each node of the triangulation of the solid. This process is transparent to the user and is performed by the class |MappingQEulerian<dim>|. In our code, we pass an object of this class as an argument to all the standard |deal.II| classes which are involved in computing the finite element values and their gradients on the deformed cells of the triangulation of the solid. In the following code snippet we illustrate this process that takes place at the beginning of the computation of the residual and of the Jacobian:

{\footnotesize%
\begin{verbatim}
...
MappingQEulerian<dim> * mapping;
...

template <int dim>
void
ImmersedFEM<dim>::residual_and_or_Jacobian(...) 
{
  if(mapping != NULL) delete mapping;

  if(par.semi_implicit == true)
      mapping = new MappingQEulerian<dim, Vector<double>, dim> 
        (par.degree, previous_xi.block(1), dh_s);
  else
    mapping = new MappingQEulerian<dim, Vector<double>, dim> 
        (par.degree, xi.block(1), dh_s);

  ...

  FEValues<dim,dim> fe_v_s_mapped (*mapping,
                                   fe_s,
                                   quad_s,
                                   update_quadrature_points);
  ...
}
\end{verbatim}}
This code snippet illustrates how to instantiate an iso-parametric mapping based on the current displacement solution, given by |xi.block(1)| or on the previous displacement solution\\ |previous_xi.block(1)|. We refer to the |deal.II| documentation of the class\\ |MappingQEulerian| for further details on the meaning of each of the arguments passed to the constructor of the class. Here it is important to notice that, once a mapping from the reference configuration to the deformed configuration is available, it is used in all instantiations of those classes which compute the values and the gradients of the basis functions on the deformed configuration (i.e., |FEValues<dim,dim>|).

Setting the parameter ``|Semi-implicit scheme|'' to |true| in the parameter file (see Section~\ref{sec:param-input-files}) will set the variable |par.semi_implicit| to |true| in the above snippet of code. The consequence of this choice is that, while the elastic response of the solid is computed at its current configuration, i.e., the Piola-Kirchhoff stress is still computed using |xi.block(1)|, the body force corresponding to this stress is applied to the fluid surrounding the body at the location |previous_xi.block(1)|, instead of |xi.block(1)|. In other words, the operator defined in Eq.~\eqref{eq: pseudo stiffness BOmega def}, and later split in the change of variable operator and in the Lagrangian elastic operator in Theorem~\ref{th: eulerian vs lagrangian stiffness} (see Eq.~\eqref{eq: EulerianLagrangianElasticity}), will use |xi.block(1)| in place of the variable $\bv w$ and |previous_xi.block(1)| in place of the variable $\bv h$. 

This splitting preserves the consistency of the method, and removes the nonlinearity due to the change of variable from the system at the cost of introducing a CFL condition on the time stepping scheme (for a more detailed discussion on this topic see \citealp{Heltai-2008-a,BoffiGastaldiHeltai-2007-a}), which ceases to be asymptotically stable. 

\subsubsection{Sparsity pattern}
\label{sec:sparsity-pattern}

A |SparsityPattern| is a |deal.II| object which stores the nonzero entries of a sparse matrix. Since we are using a |BlockSparseMatrix<double>| class to store the Jacobian of the DAE system, we need a |SparsityPattern| for each of the sub-blocks of this block. The snippet of code that generates the coupling sparsity pattern is given by

{\footnotesize%
\begin{verbatim}
  FEFieldFunction<dim, DoFHandler<dim>, Vector<double> > 
        up_field (dh_f, tmp_vec_n_dofs_up);

  vector< typename DoFHandler<dim>::active_cell_iterator > cells_f;
  vector< vector< Point< dim > > > qpoints_f;
  vector< vector< unsigned int> > maps;
  vector< unsigned int > dofs_f(fe_f.dofs_per_cell);
  vector< unsigned int > dofs_s(fe_s.dofs_per_cell);

  typename DoFHandler<dim,dim>::active_cell_iterator
    cell_s = dh_s.begin_active(),
    endc_s = dh_s.end();

  FEValues<dim,dim> fe_v_s(immersed_mapping, fe_s, quad_s,
                         update_quadrature_points);

  CompressedSimpleSparsityPattern sp1(n_dofs_up, n_dofs_W);
  CompressedSimpleSparsityPattern sp2(n_dofs_W , n_dofs_up);

  for(; cell_s != endc_s; ++cell_s)
    {
      fe_v_s.reinit(cell_s);
      cell_s->get_dof_indices(dofs_s);
      vector< Point< dim > >  &qpoints_s 
                =  fe_v.get_quadrature_points();

      up_field.compute_point_locations (qpoints_s,
                                        cells_f, qpoints_f, maps);
      for(unsigned int c=0; c<cells_f.size(); ++c)
           {
             cells_f[c]->get_dof_indices(dofs_f);
             for(unsigned int i=0; i<dofs_f.size(); ++i)
               for(unsigned int j=0; j<dofs_s.size(); ++j)
                 {
                   sp1.add(dofs_f[i],dofs_s[j]);
                   sp2.add(dofs_s[j],dofs_f[i]);
                 }
           }
    }
\end{verbatim}}
Here an |FEFieldFunction<dim>| object is constructed with a dummy finite element vector field (|tmp_vec_n_dofs_up|) to have access to its member function \\|FEFieldFunction<dim>::compute_point_locations|. This member function takes as input the location of the quadrature points in each solid cell |qpoints_s| (computed with the \\|FEValues<dim>| object |fe_v_s|, initialized with the mapping described in Section~\ref{sec:immersed-map}) and fills up a series of vectors, which allow the computation of the integrals as explained in Section~\ref{subsection: FEM implementation}. 

These vectors are respectively:
\begin{itemize}
\item |cells_f|: the vector of all \emph{fluid cells} containing at least one of the quadrature points of the immersed domain;
\item |qpoints_f|: a vector of the same length as |cells_f|, containing the custom vector of quadrature points in the \emph{fluid} reference (unit) cell, which gets transformed via the \emph{fluid mapping} to the subset of solid quadrature points |qpoints_s| (that happen to be in the cell in question); 
\item |maps|: a vector of the same length as |cells| and |qpoints_f|, which contains vectors of indices of the solid quadrature points to which the fluid quadrature points refer to, i.e., \\
  |qpoints_f[i][j]| is mapped by the \emph{fluid mapping} to the same physical location to which the point  |qpoints_s[maps[i][j]]| is mapped by the \emph{solid mapping}.
\end{itemize}

In the construction of the sparsity patterns, only the first vector, |cells_f|, is used since we only need to know which degrees of freedom are coupled. In particular, all degrees of freedom in the fluid cells contained in |cells_f| will couple with the solid cell identified with the cell iterator |cell_s|. These couplings are computed in the innermost for-loop.

\subsubsection{Residual and Jacobian}
\label{sec:residual-jacobian}

Similarly to what happens for the computation of the sparsity pattern, we use an object of type |FEFieldFunction<dim>| to compute the location of the quadrature points of the immersed solid within the fluid cells. Assembly of the coupling matrices is then possible by looping over all solid cells, and constructing custom quadrature formulas to use with the fluid cells in order to compute the integrals explained in Section~\ref{subsection: FEM implementation}. The following snippet of code explains the most relevant points:

{\footnotesize%
\begin{verbatim}
  // Loop over solid cells
  for(cell_s = dh_s.begin_active(); cell_s != endc_s; ++cell_s)
  {
    fe_v_s_mapped.reinit(cell_s);
...
    up_field.compute_point_locations (fe_v_s_mapped.get_quadrature_points(),
                                      fluid_cells,
                                      fluid_qpoints,
                                      fluid_maps);
...

// Cycle over all of the fluid cells that happen to contain some of
// the the quadrature points of the current solid cell.
      for(unsigned int c=0; c<fluid_cells.size(); ++c)
        {
          fluid_cells[c]->get_dof_indices (dofs_f);

// Local FEValues of the fluid
          Quadrature<dim> local_quad (fluid_qpoints[c]);
          FEValues<dim> local_fe_f_v (fe_f,
                                      local_quad,
                                      update_values |
                                      update_gradients |
                                      update_hessians);
          local_fe_f_v.reinit(fluid_cells[c]);
...

// Use the local_fe_f_v as you would normally do:
          for(unsigned int i=0; i<fe_s.dofs_per_cell; ++i)
            {
              unsigned int wi = i + fe_f.dofs_per_cell;
              comp_i = fe_s.system_to_component_index(i).first;
              for(unsigned int q=0; q<local_quad.size(); ++q)
                {
                  unsigned int &qs = fluid_maps[c][q];

...

                  local_res[wi] -= par.Phi_B
                                   * local_up[q](comp_i)
                                   * fe_v_s.shape_value(i,qs)
                                   * fe_v_s.JxW(qs);
...
\end{verbatim}}
In the snippet above, we show how the term $- \int_K \bv u( \bv s + \bv w(\bv s,t),t)\cdot \bv y(\bv s) d\bv s$ is assembled in practice. The point locations are computed by |up_field.compute_point_locations|. We loop over the filled vectors to compute the coupling between each of the fluid cells, |fluid_cells[c]|,  and the solid cell  |cell_s|. Since the computed quadrature points in the fluid reference cells are not standard (i.e., they are not located at Gauss quadrature points), we need to create a custom quadrature formula containing the points of interests (the object |local_quad|, initialized with |fluid_qpoints[c]|) as well as an |FEValues| object, |local_fe_f_v|, to calculate values and gradients of the fluid shape functions at the solid quadrature points.

These custom |FEValues| are then initialized with the fluid cell |fluid_cells[c]|. Notice that the correspondence between the indexing in the solid quadrature points and in the fluid custom quadrature is given by |fluid_maps[c][q]|. The rest follows the standard usage of the |deal.II| library, as can be found in any of the |deal.II| example programs.

\section{Numerics}
\label{sec:numerics}
We present in this section two numerical experiments that highlight the aspects of the accuracy and error convergence properties as well as the volume conservation feature of our numerical method.

\subsection{Static equilibrium of an annular solid comprising circumferential fibers and immersed in a stationary fluid}

This numerical test is motivated by the ones presented in \cite{BoffiGastaldiHeltaiPeskin-2008-a,GriffithLuo-2012-a}. The objective of this test is to compute the equilibrium state of an initially undeformed thick annular cylinder submerged in a stationary incompressible fluid that is contained in a rigid prismatic box having a square cross-section. Our simulation is two-dimensional and comprises an annular solid with inner radius $R$ and thickness $w$, and filled with a stationary fluid that is contained in a square box of edge length $l$ (see Fig.~\ref{fig:RingEqm-Geometry}).  In this setting, the reference and the deformed configurations of the annular solid can be conveniently described using the polar coordinate systems, whose origins coincide with the center of the annulus and whose unit vectors are given by $\left(\hat{\bv{u}}_{R}, \hat{\bv{u}}_{\Theta} \right)$ and $\left(\hat{\bv{u}}_{r}, \hat{\bv{u}}_{\theta} \right)$, respectively.
\begin{figure}[htbp]
	\begin{center}
		\includegraphics{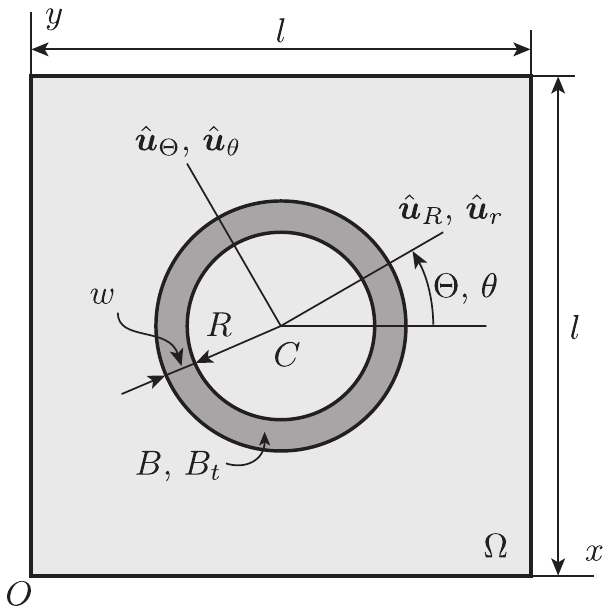}
	\caption{The reference and deformed configurations of a ring immersed in a square box filled with stationary fluid}
	\label{fig:RingEqm-Geometry}
	\end{center}
\end{figure}
This ring is located coaxially with respect to that of the box and it is subjected to the hydrostatic pressure of the fluid $p_{i}$ and $p_{o}$ at its inner and outer walls, respectively. Negligible body forces act on the system and there is no inflow or outflow of fluid across the walls of the box. Since both the solid and the fluid are incompressible, it is expected that neither the annulus nor the fluid will move at all. Therefore, the problem reduces to determining the equilibrium solution for the Lagrange multiplier field $p$. The elastic behavior of the ring is governed by a continuous distribution of concentric fibers lying in the circumferential direction. The first Piola-Kirchhoff stress for the ring is then given by
\begin{equation}
\label{eq: PK stress for ring}
\hat{\tensor{P}} = -p_{\s} \tensor{F}^{-T} + \mu^{e} \tensor{F} \hat{\bv{u}}_{\Theta} \otimes \hat{\bv{u}}_{\Theta},
\end{equation}
where $\mu^{e}$ is a constant and $p_{s}$ is the Lagrange multiplier that enforces incompressibility of the ring. As alluded to earlier, the reference configuration and the deformed configuration of the ring must coincide because of incompressibility, and the fact that the deformation of the ring must be axisymmetric in nature. For $\tensor{F}=\tensor{I}$ the constitutive response for the Cauchy stress can then be written as
\begin{equation}
\label{eq:Cauchy stress for ring}
\hat{\tensor{T}}_{\s} = -p_{\s} \tensor{I} + \mu^{e} \hat{\bv{u}}_{\theta} \otimes \hat{\bv{u}}_{\theta},
\end{equation}
where, for the deformation at hand, $\hat{\bv{u}}_{\theta} = \hat{\bv{u}}_{\Theta}$.  The balance of linear momentum for the ring can be obtained from Eq.~\eqref{eq: Cauchy theorem} as 
\begin{equation}
\label{eq: ring eqm}
-\grad \left(p_{\s}\right) + \mu^{e} \ldiv \left( \hat{\bv{u}}_{\theta} \otimes \hat{\bv{u}}_{\theta} \right) = \bv{0}.
\end{equation}
Noting that 
$\grad \left(p_{\s}\right)=\left(\partial p_{s}/\partial r \right) \hat{\bv{u}}_{r}+(1/r)(\partial p_{s}/\partial \theta) \hat{\bv{u}}_{\theta}$, 
and that 
$\ldiv \left( \hat{\bv{u}}_{\theta} \otimes \hat{\bv{u}}_{\theta} \right) = -(1/r)\hat{\bv{u}}_{r}$,
Eq.~\eqref{eq: ring eqm} can be rewritten as 
\begin{equation}
\label{eq: ring eqm 2}
-\frac{\partial p_{\s}}{\partial r} - \frac{\mu^{e}}{r} = 0
\qquad
\mathrm{and}
\qquad
\frac{\partial p_{\s}}{\partial \theta} = 0.
\end{equation}
From Eq.~\eqref{eq: ring eqm 2}, it can be concluded that the Lagrange multiplier enforcing incompressibility $p_{\s}$ is an axisymmetric function of the form 
\begin{equation}
\label{eq: solution for p_s}
p_{\s} = c - \mu^{e} \ln \left(\frac{r}{R} \right),
\end{equation}
where $c$ is a constant. The satisfaction of the traction boundary conditions at the inner and outer walls of the ring demand that $p_{\s}\vert_{r=R}=p_{i}$ and $p_{\s}\vert_{r=R+w}=p_{o}$ and hence we can obtain that
\begin{equation}
\label{eq: solution for p_s 2}
p_{\s} = p_{i} - \mu^{e} \ln \left(\frac{r}{R}\right),
\qquad
p_{o} = p_{i} - \mu^{e} \ln \left(1+\frac{w}{R}\right)
\end{equation}
Note that Lagrange multiplier $p$ defined over the control volume corresponds to $p_{\s}$ in the region occupied by the solid. By constraining the average value of $p$ over the entire control volume to be zero we arrive at the following solution for the equilibrium problem:
\begin{equation}
p=
\begin{cases}
p_{o}=-\frac{\pi \mu^{e}}{2 l^{2}} \left( \left(R+w\right)^{2}-R^{2}\right)  & \mathrm{for} \quad R + w \leq r,\\
p_{s}=\mu^{e} \ln (\frac{R+w}{r})-\frac{\pi \mu^{e}}{2 l^{2}} \left( \left(R+w\right)^{2}-R^{2}\right) & \mathrm{for} \quad R <r < R+w,\\
p_{i}=\mu^{e} \ln (1+\frac{w}{R})-\frac{\pi \mu^{e}}{2 l^{2}} \left( \left(R+w\right)^{2}-R^{2}\right) &\mathrm{for} \quad  r \leq R,
\end{cases}
\label{eqn: p for ring eqm}
\end{equation}
with velocity of fluid $\bv{u}=\bv{0}$ and the displacement of the solid $\bv{w}=\bv{0}$. Note that Eq.~\eqref{eqn: p for ring eqm} is different from Eq.~(69) of \cite{BoffiGastaldiHeltaiPeskin-2008-a}, where $p$ varies linearly with $r$ (we believe this to be in error).

For all our numerical simulations we have used $R =\np[m]{0.25}$, $w=\np[m]{0.06250}$, $l=\np[m]{1.0}$ and $\mu^{e}=\np[Pa]{1}$ and for these values we obtain $p_{i}=\np[Pa]{0.16792}$ and $p_{o}=\np[Pa]{-0.05522}$ using Eq.~\eqref{eqn: p for ring eqm}. We have used $\rho=\np[kg/m^{3}]{1.0}$, dynamic viscosity $\mu=\np[Pa \!\cdot\! s]{1.0}$, and time step size $h=\np[s]{1e-3}$ in our tests. For all our numerical tests we have used $Q2$ elements to represent $\bv{w}$ of the solid, whereas we have used (i) $Q2/P1$ elements, and (ii) $Q2/Q1$ elements to represent $\bv{v}$ and $p$ over the control volume. We present a sample profile of $p$ over the entire control volume and its variation along different values of $y$, after one time step, in Fig.~\ref{fig:DGP_Pressure} and Fig.~\ref{fig:FEQ_Pressure} for $Q2/P1$ and $Q2/Q1$ elements, respectively.

We assess the convergence property of our numerical scheme by obtaining the convergence rate of the error between the exact and the numerical solutions of this equilibrium problem. The order of the rate of convergence (see, Tables~\ref{tab:ring-eqm-dgp} and \ref{tab:ring-eqm-feq} for $Q2/P1$ and $Q2/Q1$ elements, respectively) is 2.5 for the $L^{2}$ norm of the velocity, 1.5 for the $H^{1}$ norm of the velocity and 1.5 for the $L^{2}$ norm of the pressure which matches the rates presented in \cite{BoffiGastaldiHeltaiPeskin-2008-a}. In all these numerical tests we have used 1856 cells with 15776 DoFs for the solid.

The parameter files used for these tests can be found under the directory\\|prms/RingEqm_XXX_fref_Y_param.prm|, where |XXX| is either |dgp| or |feq| and |Y| is 4, 5, 6 or 7, according to the type of pressure finite element and to the fluid refinement level. The tests can be run under the |step-feibm| directory, by typing \\
|./step-feibm prms/RingEqm_XXX_fref_Y_param.prm|

\begin{figure}[htbp]
	\begin{center}
	\subfigure[Over the entire domain]
	{\includegraphics[width=3in]{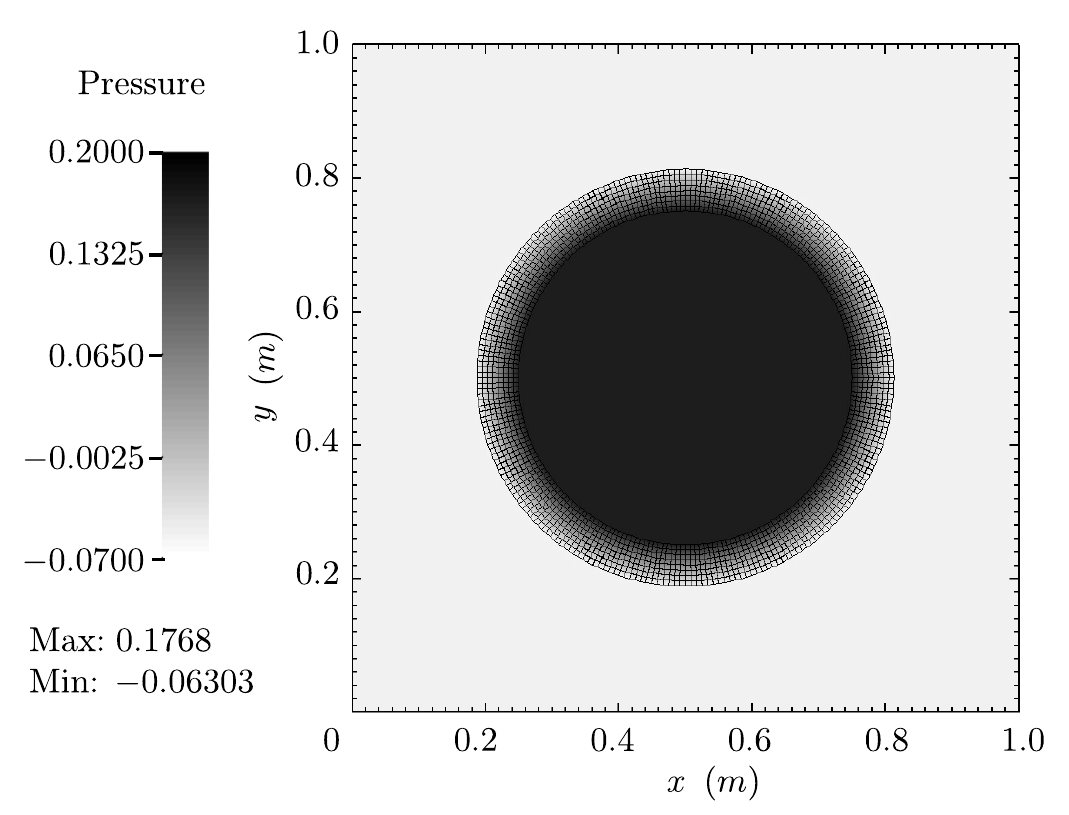}
	}
	\subfigure[At different values of $y$]
	{\includegraphics[width=3.1in]{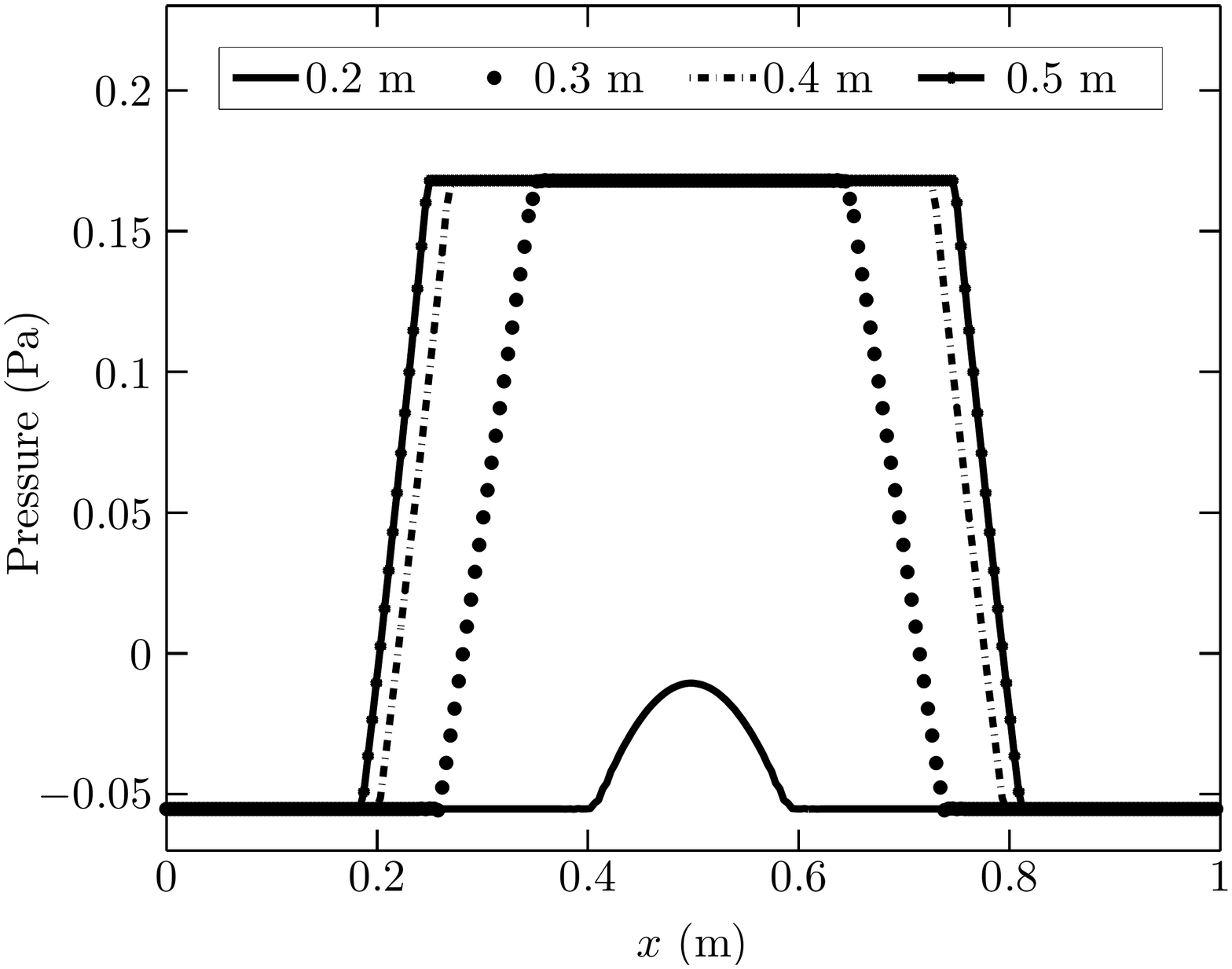}
	}
	\caption{The values of $p$ after one time step when using $P1$ elements for $p$}
	\label{fig:DGP_Pressure}
	\end{center}
\end{figure}
\begin{table}[htbp]\small
\caption{Error convergence rate obtained when using $P1$ element for $p$ after one time step}
\begin{center}
\begin{tabular*} {\textwidth} {@{\extracolsep{\fill}} c c c c c c c c}
\toprule
No. of cells &
No. of DoFs  &
$\|\mathbf{u}_{h}-\mathbf{u}\|_{0}$ & &
$\|\mathbf{u}_{h}-\mathbf{u}\|_{1}$ & &
$\|p_{h}-p\|_{0}$ & 
\\
\midrule
 256 &   2946 & 2.00605e-05 & -    & 1.95854e-03 & -    & 6.71603e-03 & -    \\
1024 &  11522 & 3.69389e-06 & 2.44 & 7.44696e-04 & 1.40 & 2.47476e-03 & 1.44 \\
4096 &  45570 & 5.76710e-07 & 2.68 & 2.25134e-04 & 1.73 & 8.74728e-04 & 1.50 \\
16384& 181250 & 1.06127e-07 & 2.44 & 8.24609e-05 & 1.45 & 3.14028e-04 & 1.48 \\ 
\bottomrule
\end{tabular*}
\end{center}
\label{tab:ring-eqm-dgp}
\end{table}
\begin{figure}[htbp]
	\begin{center}
	\subfigure[Over the entire domain]
	{\includegraphics[width=3in]{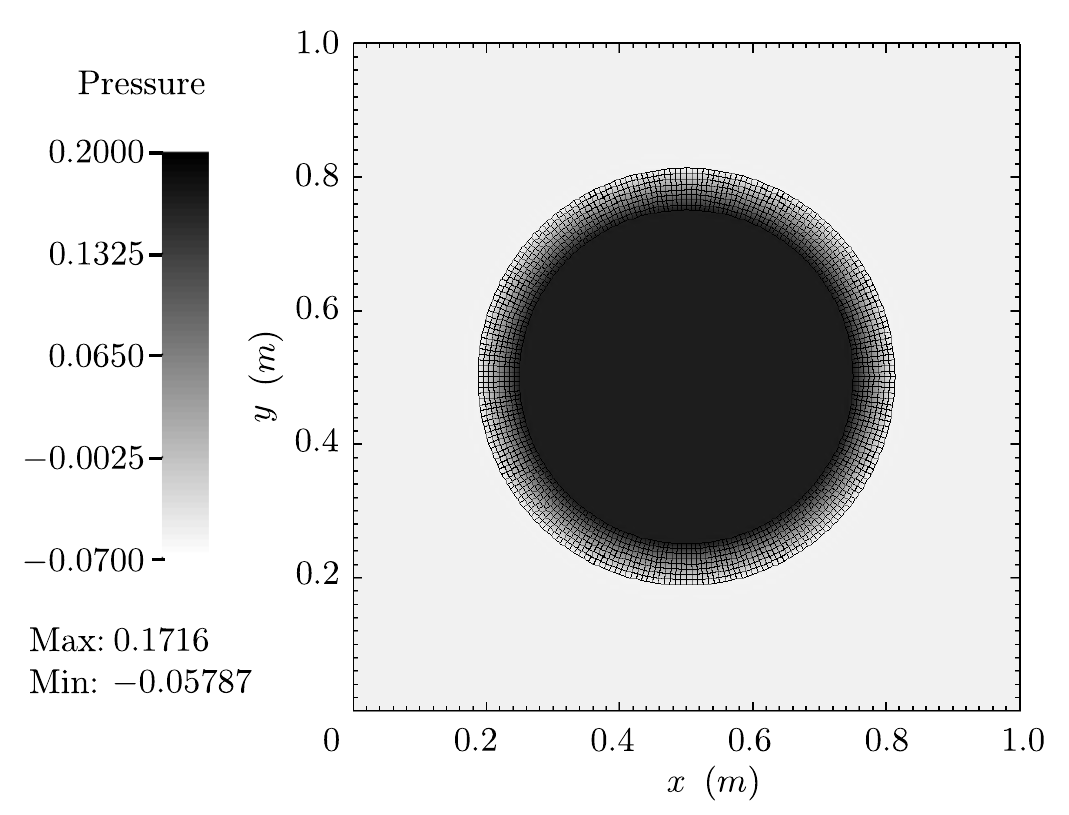}
	}
	\subfigure[At different values of $y$]
	{\includegraphics[width=3.1in]{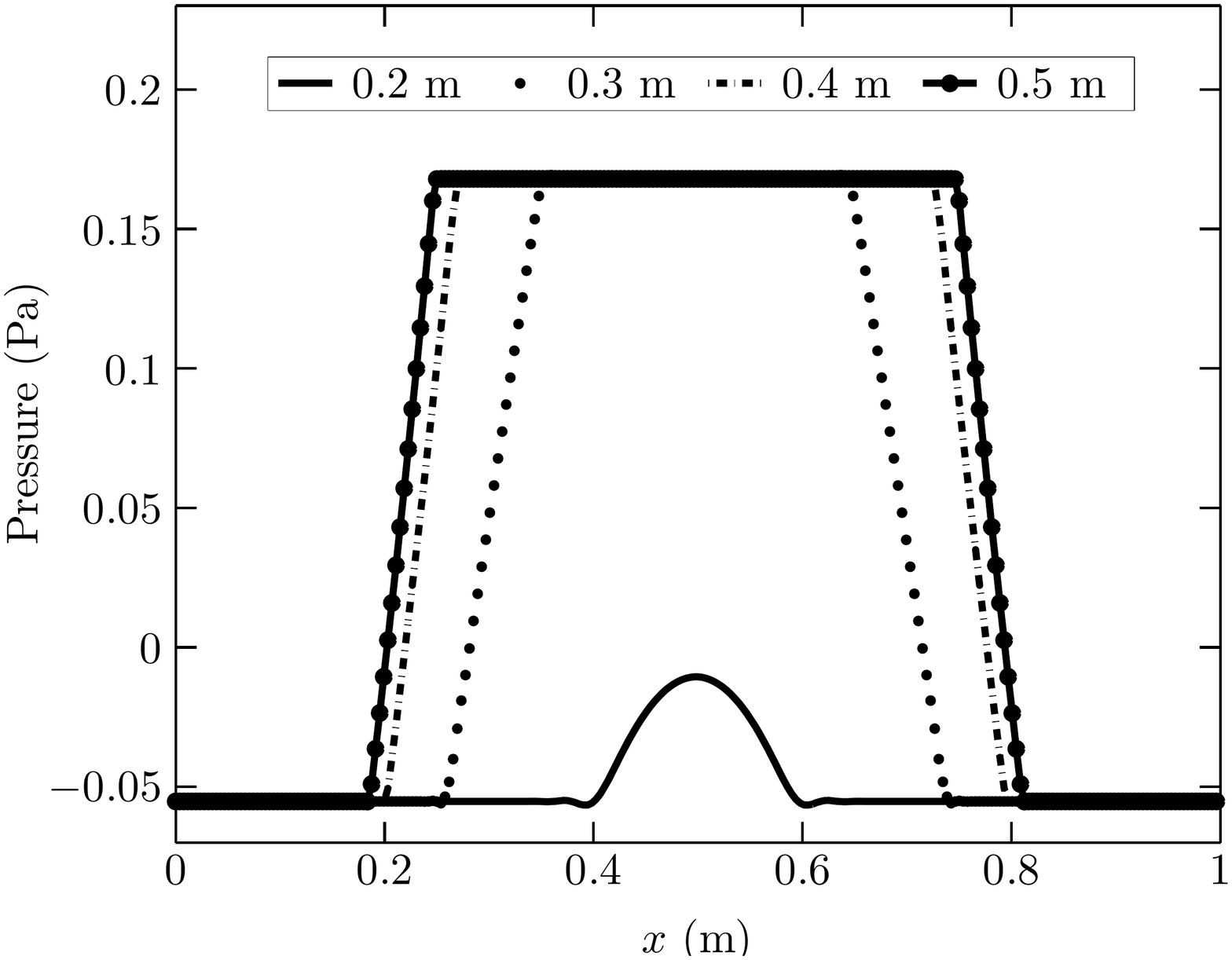}
	}
	\caption{The values of $p$ after one time step when using $Q1$ elements for $p$}
	\label{fig:FEQ_Pressure}
	\end{center}
\end{figure}
\begin{table}[htbp]\small
\caption{Error convergence rate obtained when using $Q1$ element for $p$ after one time step}
\begin{center}
\begin{tabular*} {\textwidth} {@{\extracolsep{\fill}} c c c c c c c c}
\toprule
No. of cells &
No. of DoFs  &
$\|\mathbf{u}_{h}-\mathbf{u}\|_{0}$ & &
$\|\mathbf{u}_{h}-\mathbf{u}\|_{1}$ & &
$\|p_{h}-p\|_{0}$ & 
\\
\midrule
 256 &   2467 & 4.36912e-05 &    - & 2.79237e-03 &    - & 7.39310e-03 &- \\ 
1024 &   9539 & 6.14959e-06 & 2.83 & 9.02397e-04 & 1.63 & 2.42394e-03 & 1.61 \\ 
4096 &  37507 & 1.28224e-06 & 2.26 & 3.49329e-04 & 1.37 & 9.10608e-04 & 1.41 \\ 
16384 & 148739 & 2.33819e-07 & 2.46& 1.25626e-04 & 1.48 & 3.27256e-04 & 1.48 \\
\bottomrule
\end{tabular*}
\end{center}
\label{tab:ring-eqm-feq}
\end{table}

\subsection{Disk entrained in a lid-driven cavity flow }

We test the volume conservation of our numerical method by measuring the change in the area of a disk that is entrained in a lid-driven cavity flow of an incompressible, linearly viscous fluid. This test is motivated by similar ones presented in \cite{WangZhang_2010_Interpolation-functions-0,GriffithLuo-2012-a}. Referring to Fig.~\ref{fig:LDCFlow-Ball-Geometry}, the disk has a radius $R=\np[m]{0.2}$ and its center $C$ is initially positioned at $x=\np[m]{0.6}$ and $y=\np[m]{0.5}$ in the square cavity whose each edge has the length $l=\np[m]{1.0}$. Body forces on the system are negligible. The two different constitutive models for the elastic response of the disk which we have used for our simulations are as follows:
\begin{align}
\mbox{case 1: }\quad \hat{\tensor{P}} &= -p_{\s} \tensor{I} + \mu^{e} \left(\tensor{F} -\tensor{F}^{-\mathrm{T}}\right),
\label{eq: INH0} \\
\mbox{case 2: }\quad \hat{\tensor{P}} &= -p_{\s} \tensor{I} + \mu^{e} \tensor{F}. 
\label{eq: INH1}
\end{align}
We have used the following parameters: $\rho=\np[kg/m^{3}]{1.0}$, dynamic viscosity $\mu=\np[Pa\!\cdot\! s]{0.01}$, shear modulus $\mu^{e} = \np[Pa]{0.1}$ and $U=\np[m/s]{1.0}$. For our numerical simulations we have used $Q2$ elements to represent $\bv{w}$ of the disk whereas we have used $Q2/P1$ element for the fluid. The disk is represented using 320 cells with 2626 DoFs and the control volume has 4096 cells and 45570 DoFs. The time step size $h=\np[s]{1e-2}$. We consider the time interval $0<t \leq \np[s]{8}$ during which the disk is lifted from its initial position along the left vertical wall, drawn along underneath the lid and finally dragged downwards along the right vertical wall of the cavity (see, Figs.~\ref{fig:LDCFlowBall-DGP-ZeroResStress} and \ref{fig:LDCFlowBall-DGP-ResStress}). As the disk trails beneath the lid, it experiences large shearing deformations (see, Figs.~\ref{fig:LDCFlowBall-DGP-ZeroResStress-Deformation} and \ref{fig:LDCFlowBall-DGP-ResStress-Deformation}). Ideally the disk should have retained its original area over the course of time because the incompressibility of the media and the nature of the motion require that the disk change its shape only and not its volume. However, from our numerical scheme we obtain an area change of the disk of about $6\%$ for case 1 (see Fig.~\ref{fig:LDCFlowBall-DGP-ZeroResStress-AreaChange}) and about $4\%$ for case 2 (see Fig.~\ref{fig:LDCFlowBall-DGP-ResStress-AreaChange}).

The parameter files used for these two tests can be found under the directory \\|deal.II/step-feibm/prms|, and are named |LDCFlow_Ball_DGP_INH0_param.prm| and \\|LDCFlow_Ball_DGP_INH1_param.prm| respectively. The tests can be run under the \\
|deal.II/step-feibm| directory, by typing \\
|./step-feibm prms/LDCFlow_Ball_DGP_INH0_param.prm|\\
and \\
|./step-feibm prms/LDCFlow_Ball_DGP_INH1_param.prm|

\begin{figure}[htbp]
	\begin{center}
		\includegraphics{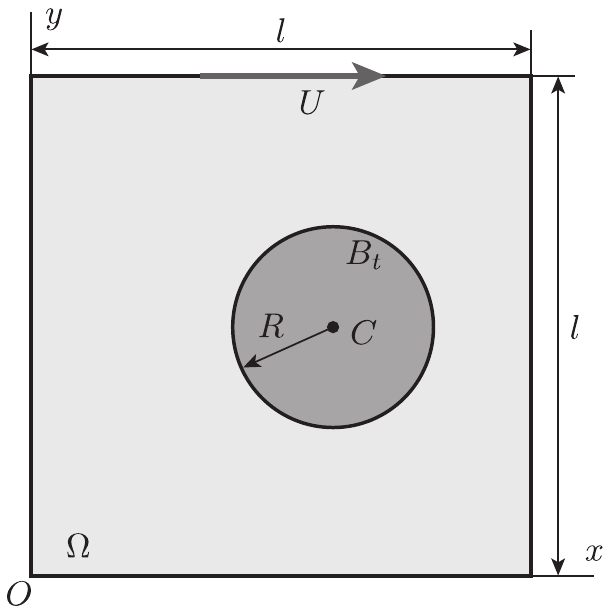}
		\caption{The initial configuration of an immersed disk entrained in a flow in a square cavity whose lid is driven with a velocity $U$ towards the right}
	\label{fig:LDCFlow-Ball-Geometry}
	\end{center}
\end{figure}
\begin{figure}[htbp]
	\begin{center}
	\subfigure
	{\includegraphics[trim=15 10 10 20, clip=true]{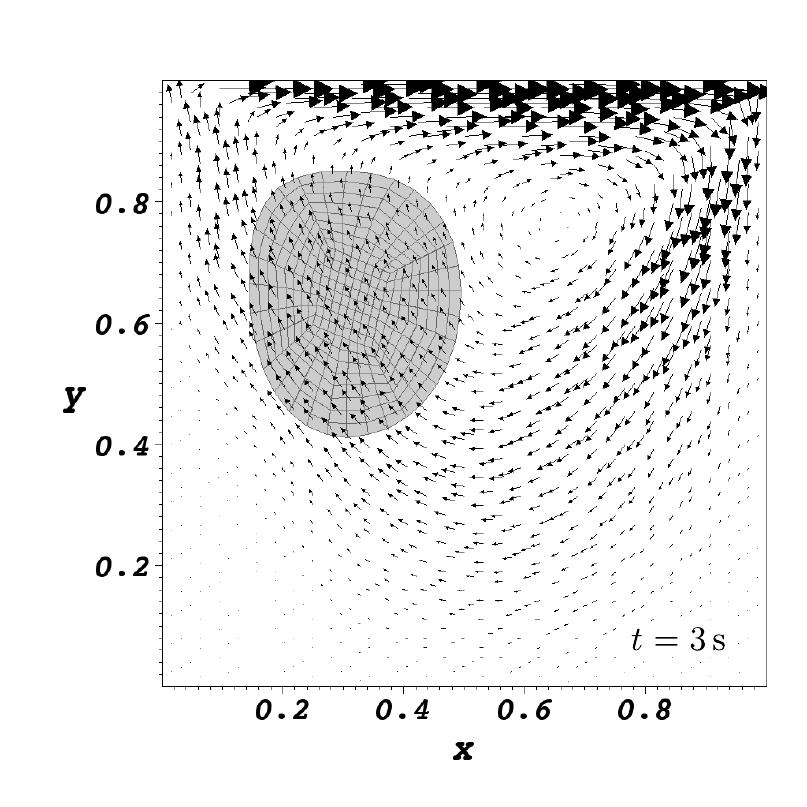}
	\label{fig:INH0-velocity-t3s}
	}
	\subfigure
	{\includegraphics[trim=15 10 10 20, clip=true]{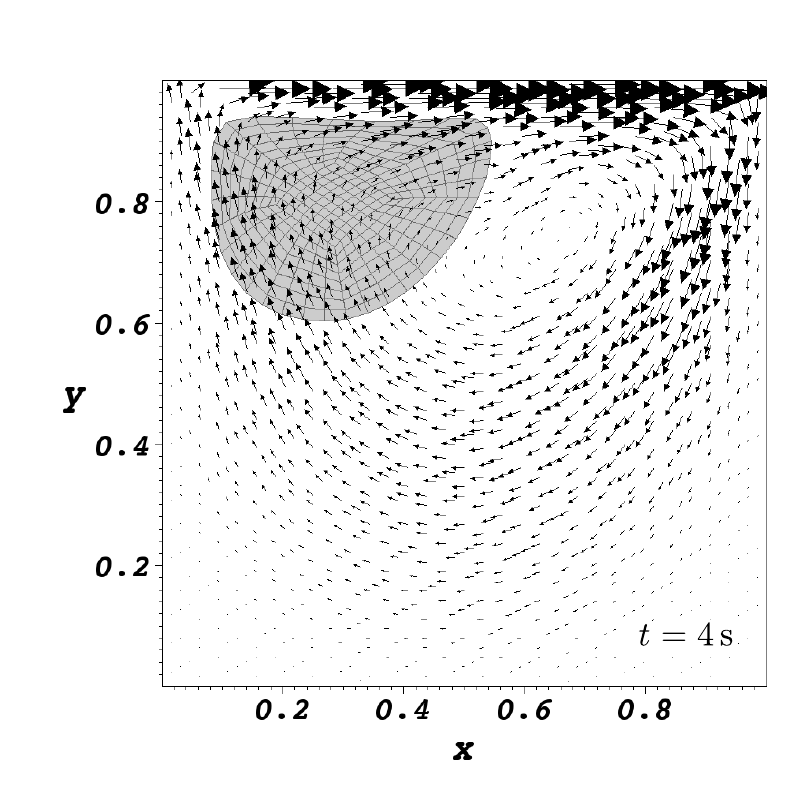}
	\label{fig:INH0-velocity-t4s}
	} 
	\subfigure
	{\includegraphics[trim=15 10 10 20, clip=true]{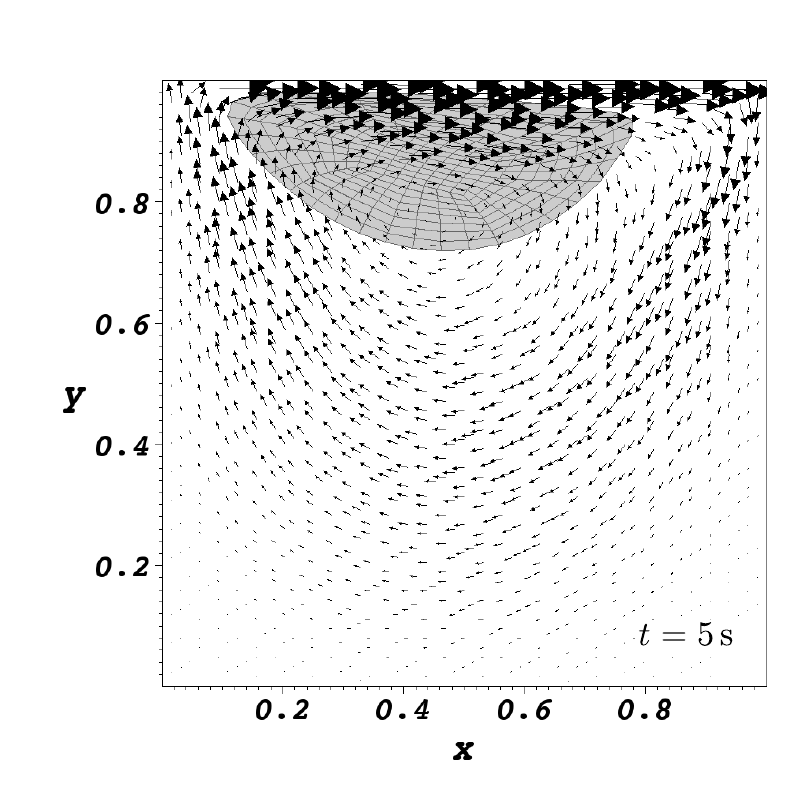}
	\label{fig:INH0-velocity-t5s}
	}
	\subfigure
	{\includegraphics[trim=15 10 10 20, clip=true]{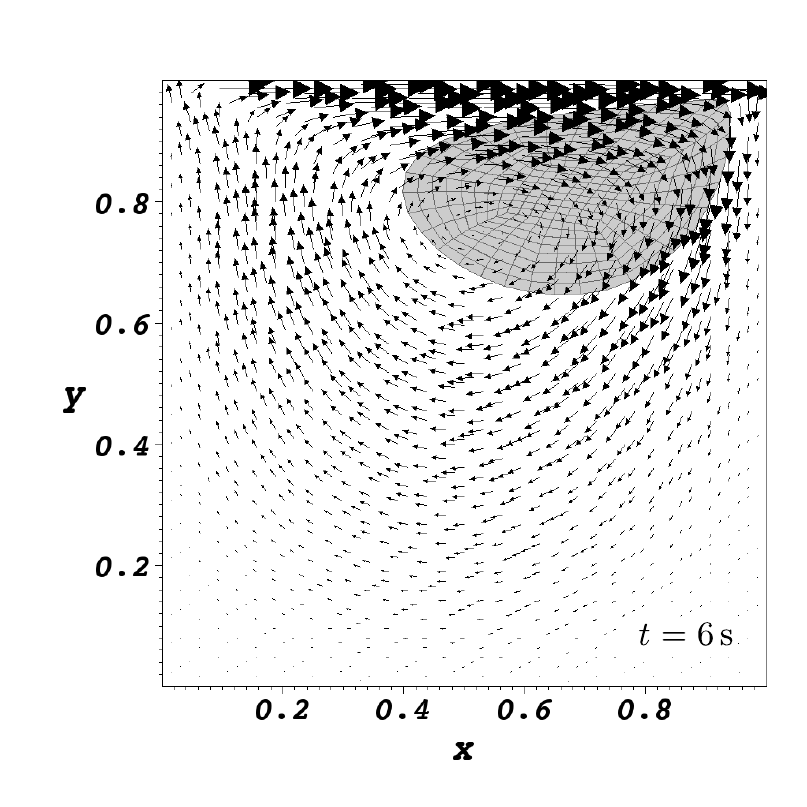}
	\label{fig:INH0-velocity-t6s}
	}
	\subfigure{
	\includegraphics[trim=15 10 10 20, clip=true]{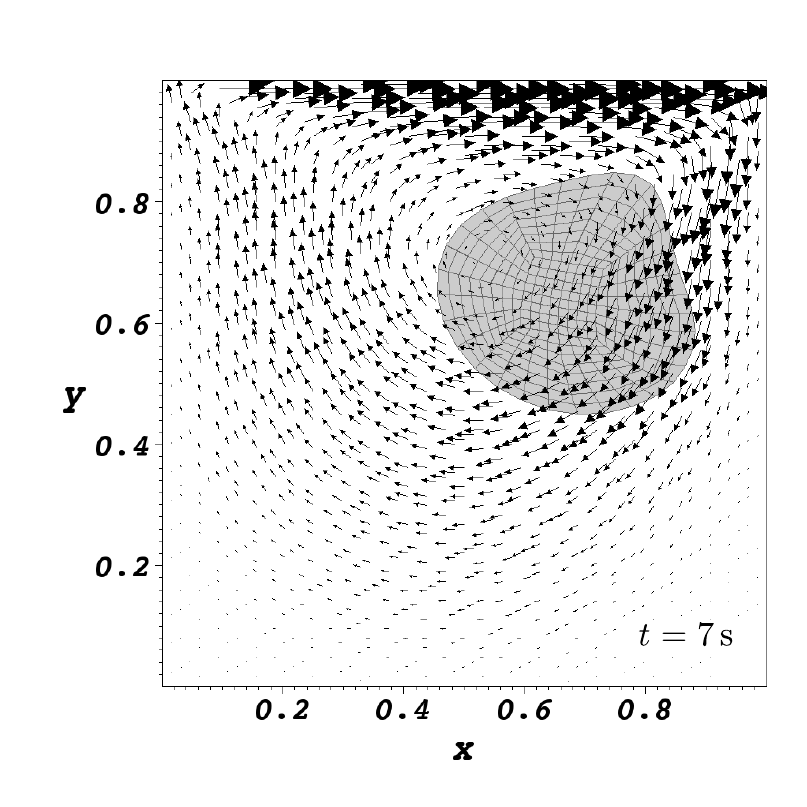}
	\label{fig:INH0-velocity-t7s}
	} 
	\subfigure{
	\includegraphics[trim=15 10 10 20, clip=true]{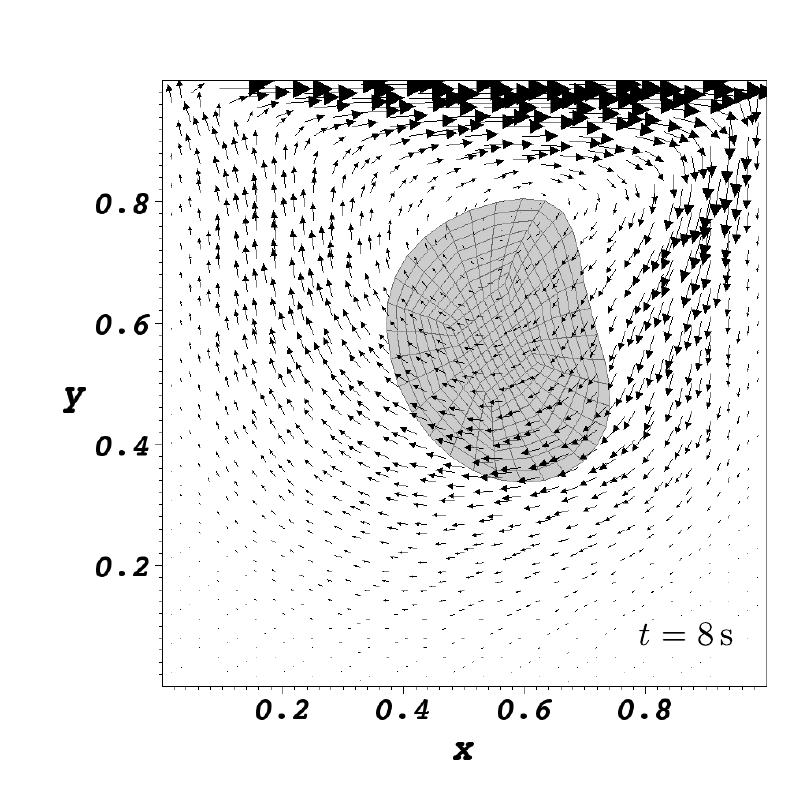}
	\label{fig:INH0-velocity-t8s}
	} 
\caption{The motion of the disk for case 1 at different instants of time}
\label{fig:LDCFlowBall-DGP-ZeroResStress}
\end{center}
\end{figure}
\begin{figure}[htbp]
	\begin{center}
	\subfigure
	{\includegraphics[trim=0 10 0 70, clip=true]{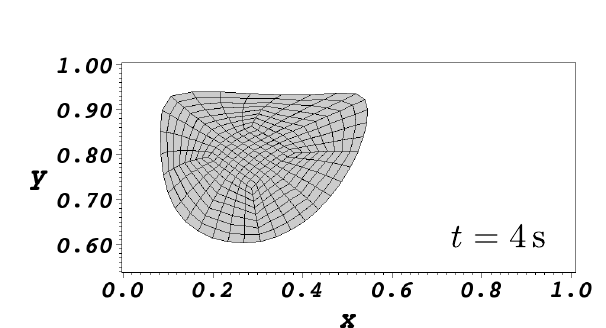} 
	\label{fig:INH0-solidlocation-t4s}
	}
	\subfigure
	{\includegraphics[trim=0 10 0 70, clip=true]{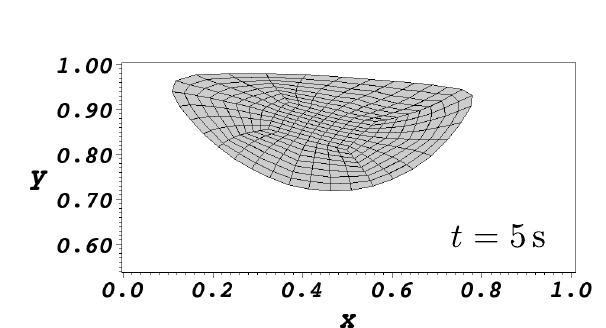}
	\label{fig:fig:INH0-solidlocation-t5s}
	} 
	\subfigure
	{\includegraphics[trim=0 10 0 70, clip=true]{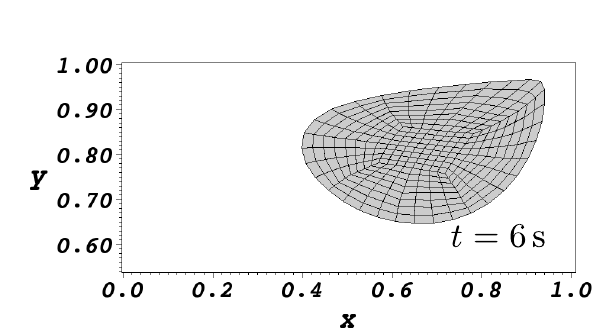}
	\label{fig:fig:INH0-solidlocation-t6s}
	}
	\caption{Enlarged view of the disk for case 1 depicting its shape and location at various instants of time}
	\label{fig:LDCFlowBall-DGP-ZeroResStress-Deformation}
	\end{center}
\end{figure}
\begin{figure}[htbp]
\begin{center}
\includegraphics[scale=0.5, trim=0.2in 2.5in 0.5in 2.5in, clip=true]{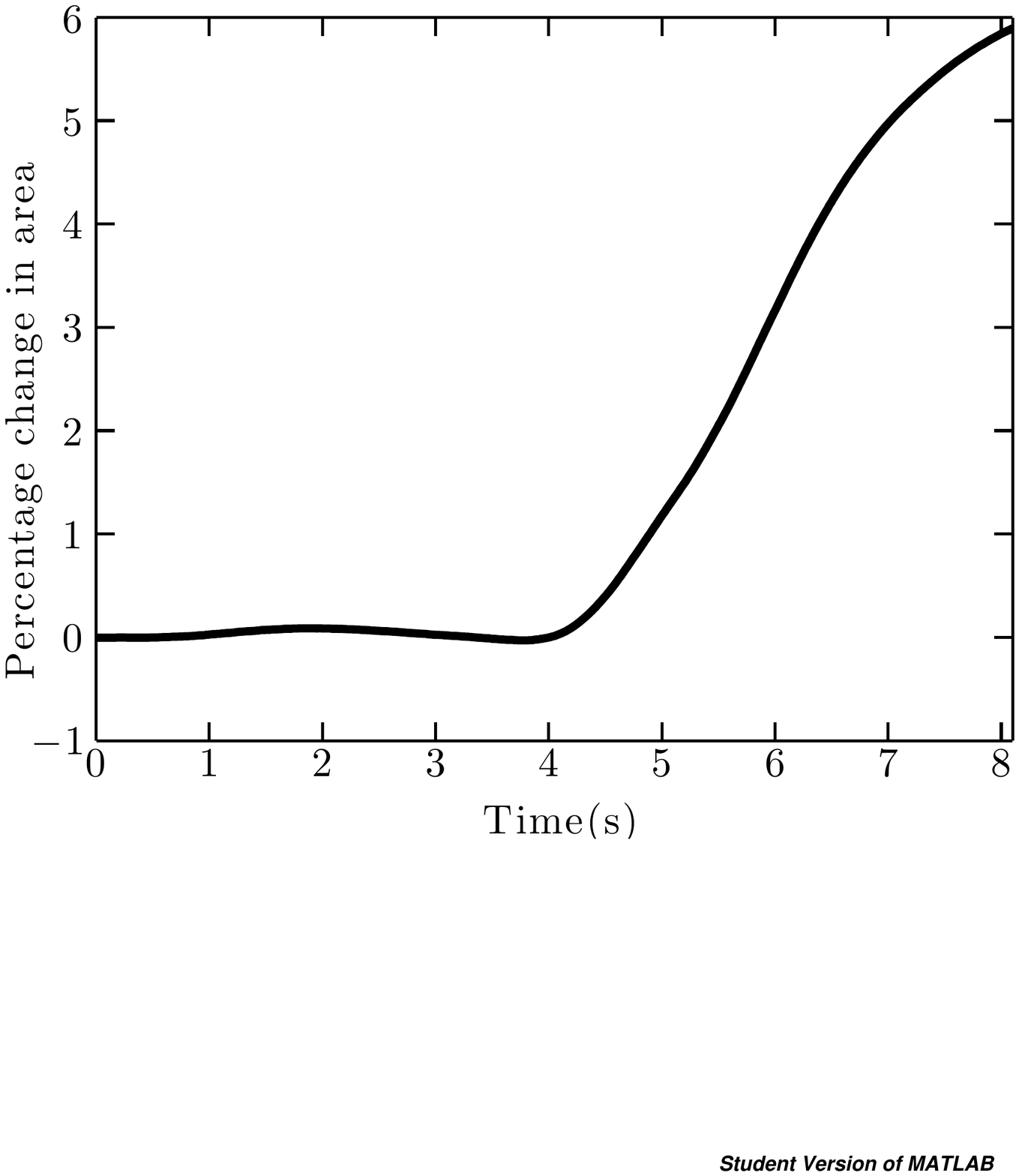}
\caption{The percentage change in the area of the disk for case 1 over time}
\label{fig:LDCFlowBall-DGP-ZeroResStress-AreaChange}
\end{center}
\end{figure}
\begin{figure}[htbp]
	\begin{center}
	\subfigure
	{\includegraphics[trim=15 10 10 20, clip=true]{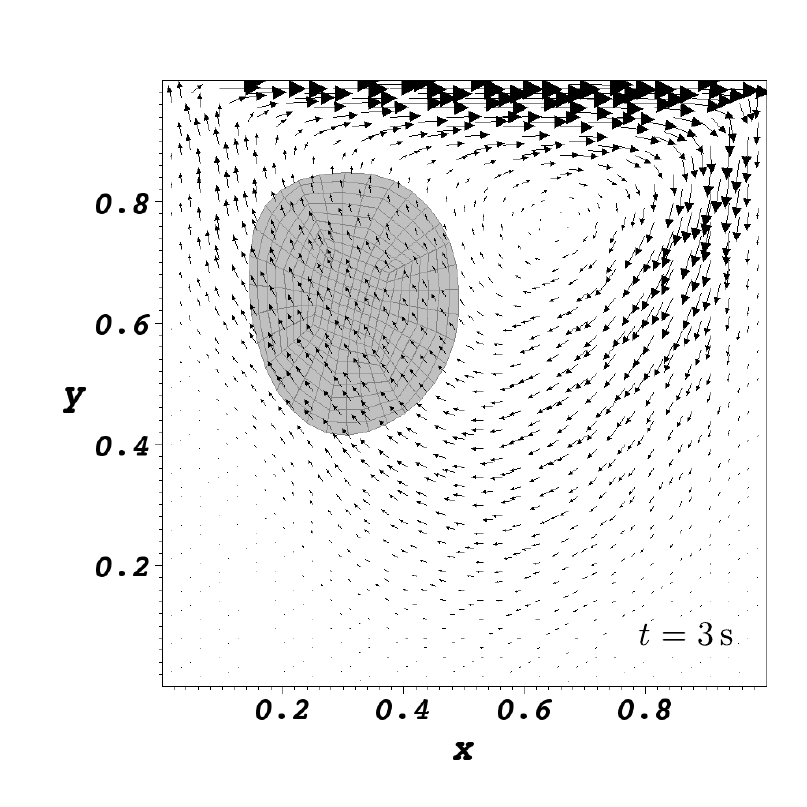}
	\label{fig:INH1-velocity-t3s}
	}
	\subfigure
	{\includegraphics[trim=15 10 10 20, clip=true]{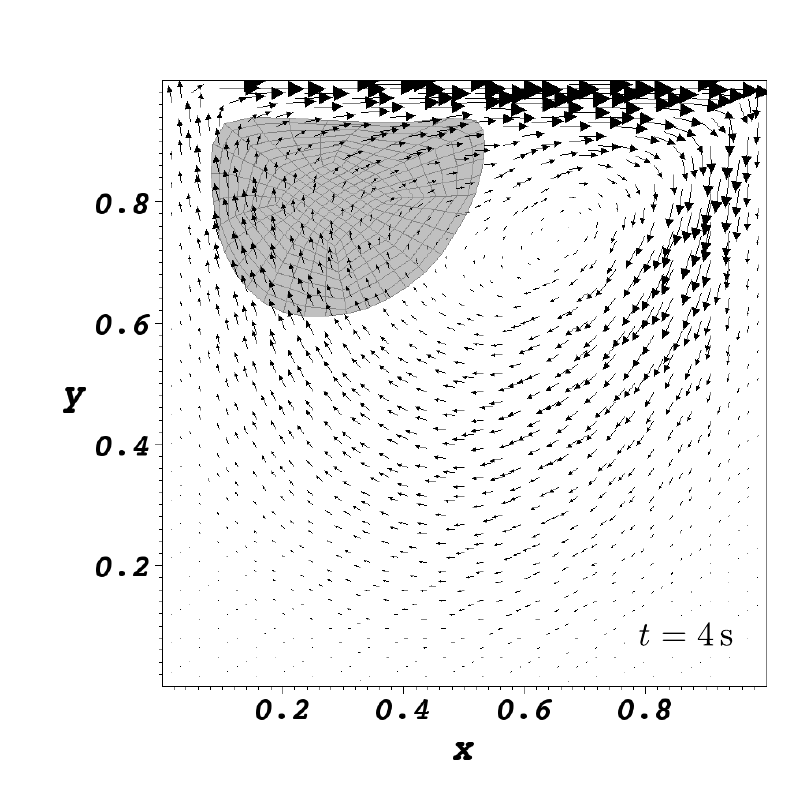}
	\label{fig:INH1-velocity-t4s}
	} 
	\subfigure
	{\includegraphics[trim=15 10 10 20, clip=true]{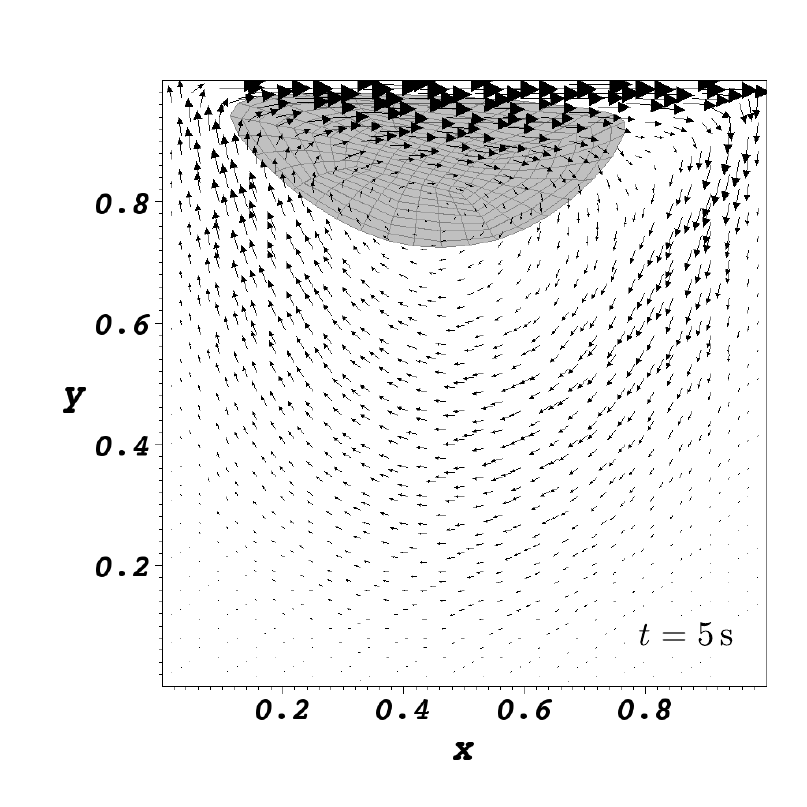}
	\label{fig:INH1-velocity-t5s}
	}
	\subfigure
	{\includegraphics[trim=15 10 10 20, clip=true]{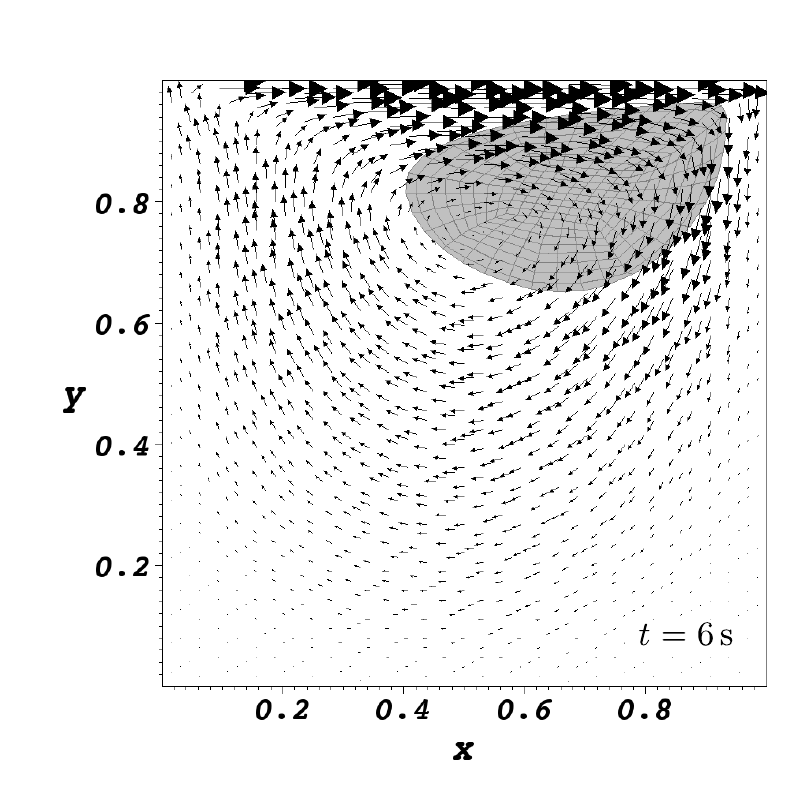}
	\label{fig:INH1-velocity-t6s}
	}
	\subfigure{
	\includegraphics[trim=15 10 10 20, clip=true]{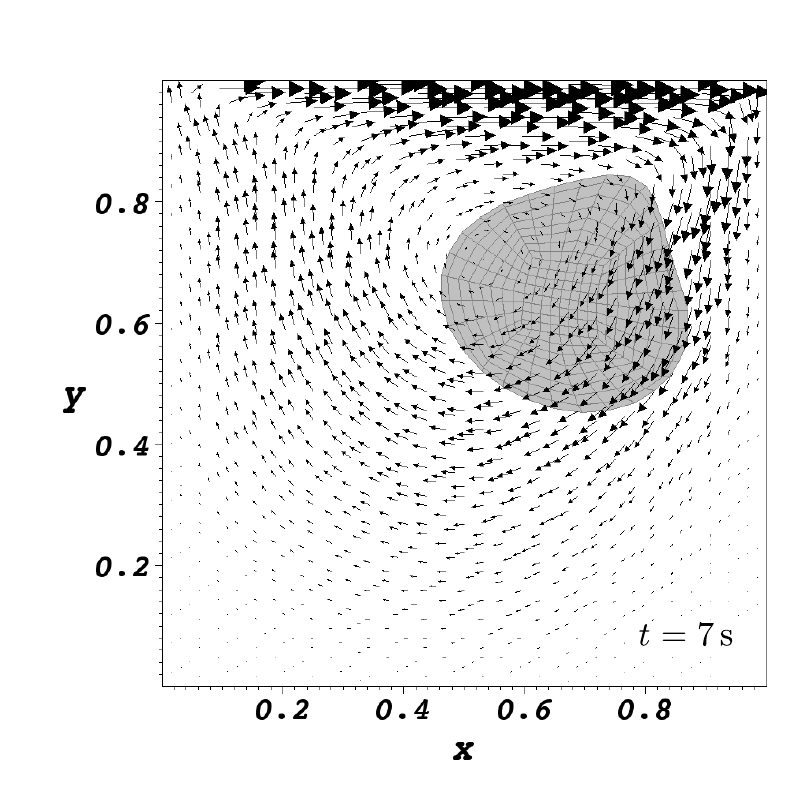}
	\label{fig:INH1-velocity-t7s}
	} 
	\subfigure{
	\includegraphics[trim=15 10 10 20, clip=true]{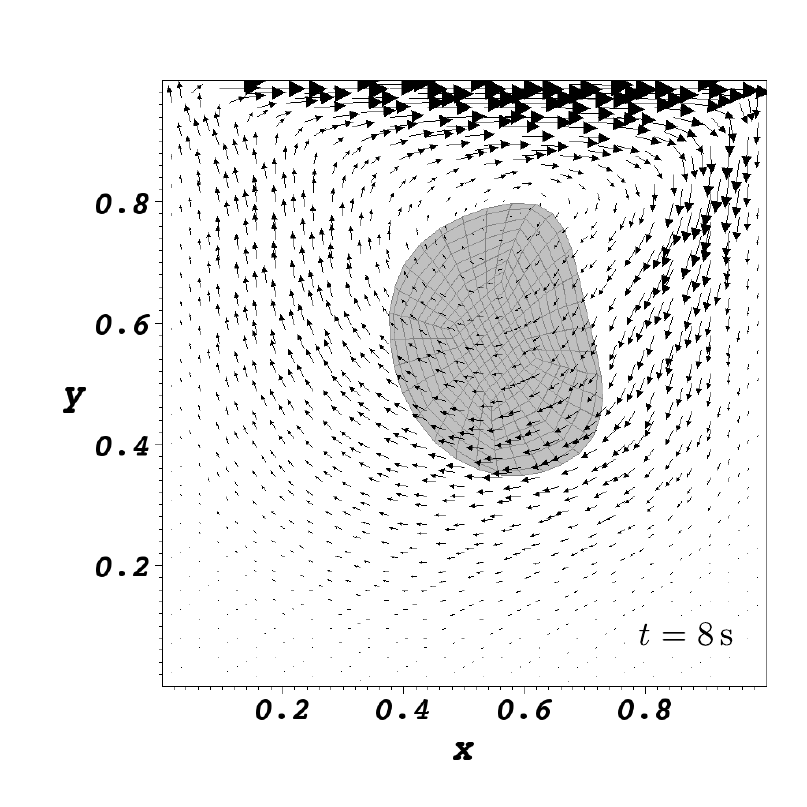}
	\label{fig:INH1-velocity-t8s}
	} 
\caption{The motion of a disk for case 2 at different instants of time}
\label{fig:LDCFlowBall-DGP-ResStress}
\end{center}
\end{figure}
\begin{figure}[htbp]
	\begin{center}
	\subfigure
	{\includegraphics[trim=0 10 0 70, clip=true]{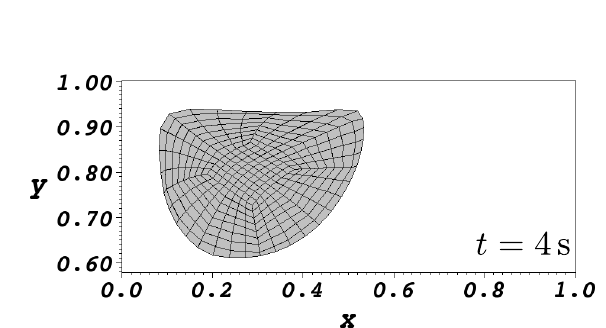} 
	\label{fig:INH1-solidlocation-t4s}
	}
	\subfigure
	{\includegraphics[trim=0 10 0 70, clip=true]{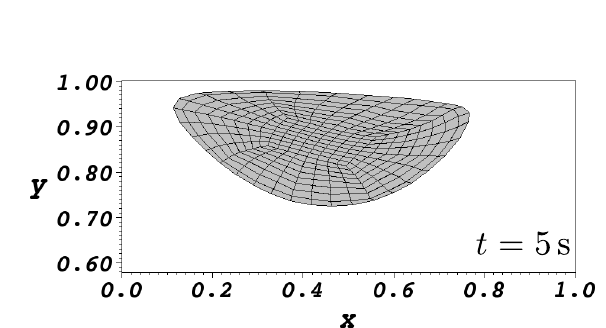}
	\label{fig:fig:INH1-solidlocation-t5s}
	} 
	\subfigure
	{\includegraphics[trim=0 10 0 70, clip=true]{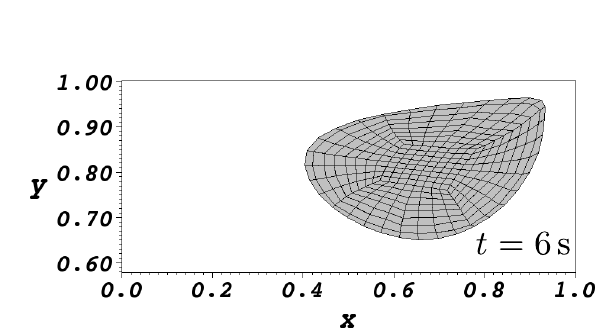}
	\label{fig:fig:INH1-solidlocation-t6s}
	}
	\caption{Enlarged view of the disk for case 2 depicting its shape and location at various instants of time}
	\label{fig:LDCFlowBall-DGP-ResStress-Deformation}
	\end{center}
\end{figure}
\begin{figure}[htbp]
\begin{center}
\includegraphics[scale=0.5, trim=0.2in 2.5in 0.5in 2.5in, clip=true]{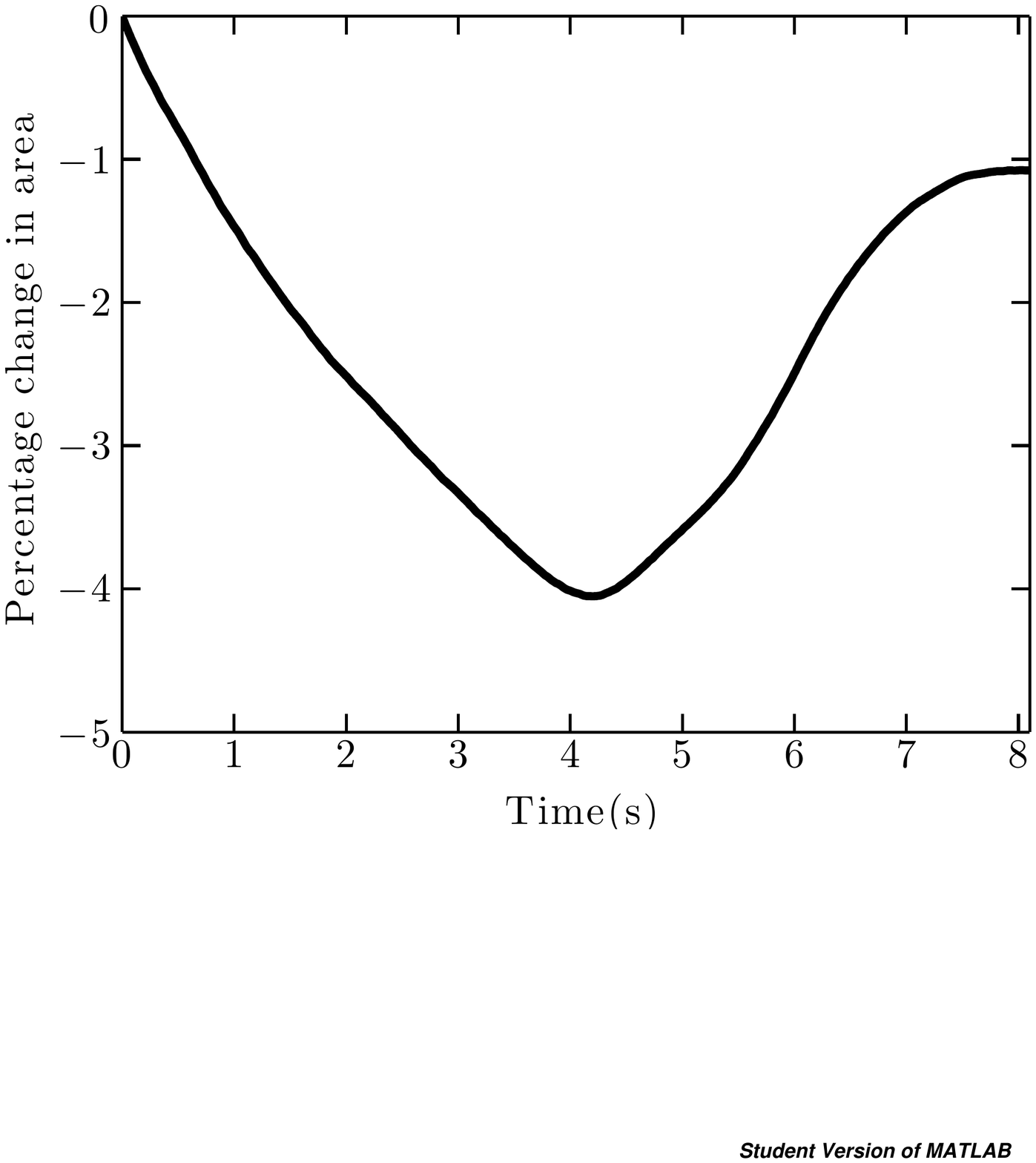}
\caption{The percentage change in the area of the disk for case 2 over time}
\label{fig:LDCFlowBall-DGP-ResStress-AreaChange}
\end{center}
\end{figure}


\section*{Acknowledgements}
\label{sec:acknowledgements}

This work was partially supported by the  ``Young Scientist Grant'' number ANFU.685, made available by the ``Scuola Internazionale Superiore di Studi Avanzati'' to the first author.

}

\bibliographystyle{FG-AY-bibstyle}                  %
\bibliography{ch_feibm_ANS} %

\end{document}